\documentclass[12pt]{amsart}
\usepackage{amsfonts,amsmath,amssymb}
\numberwithin{equation}{section}

\newtheorem{lemma}{Lemma}[section]
\newtheorem{theorem}{Theorem}[section]
\newtheorem{cor}{Corollary}[section]
\newtheorem{prop}{Proposition}[section]
\newtheorem{remark}{Remark}[section]
\newtheorem{example}{Example}[section]

\begin{document}

\def\g{GL(n,K)}
\def\gz{GL(n,\z)}
\def\Z{{\mathbb Z}}
\def\Q{{\mathbb Q}}
\def\C{{\mathbb C}}
\def\d{\nabla(\lambda)}
\def\dz{\nabla_\z(\lambda)}
\def\l{L(\lambda)}
\def\m{M(\lambda)}
\def\ov{\overline}
\def\s{S_K(n,r)}
\def\la{\lambda}
\let\tilde=\widetilde
\def\dt{\nabla(\widetilde\lambda)}
\def\Hom{\mathop{\rm Hom}\nolimits}
\def\End{\mathop{\rm End}\nolimits}
\def\im{\mathop{\rm im}\nolimits}
\def\tr{\mathop{\rm tr}\nolimits}

\def\cI{{\mathcal I}}
\def\cJ{{\mathcal J}}
\def\cK{{\mathcal K}}
\def\cT{{\mathcal T}}
\def\cB{{\mathcal B}}
\def\cD{{\mathcal D}}
\def\cP{{\mathcal P}}
\def\cE{{\mathcal E}}
\def\cF{{\mathcal F}}
\def\cS{{\mathcal S}}
\def\cT{{\mathcal T}}

\def\bi{\overline{\imath}}
\def\bj{\overline{\jmath}}

\def\boldi{{\mathbf i}}
\def\bolde{{\mathbf e}}

\def\bm{{\bf m}}

\def\sgn{\hbox{sgn}}

\title{A basis of bideterminants for the coordinate ring of the orthogonal
group}

\author{Gerald Cliff}

\address{University of Alberta, Department of Mathematical and Statistical
Sciences, Edmonton, Alberta, Canada T6G 2G1}
\email{gcliff@math.ualberta.ca}
\thanks{This research was supported in part by a grant from the
Natural Sciences and Engineering Research Council of Canada.}

\begin{abstract}We give a basis of bideterminants for the
coordinate ring $K[O(n)]$  of the
orthogonal group $O(n,K)$, where $K$ is an infinite field of characteristic
not 2. The bideterminants are indexed by pairs of Young tableaux which
are $O(n)$-standard in the sense of King-Welsh. We also give an explicit
filtration of $K[O(n)]$ as an $O(n,K)$-bimodule, whose factors are
isormorphic to the tensor product of orthogonal analogues of left
and right Schur modules.
\end{abstract}

\maketitle

\section{Introduction}
Let $K$ be an infinite field. It was shown by Mead [M] and by
Doubilet-Rota-Stein [DRS] that the polynomial ring $K[X(i,j)]$
in $n^2$ variables has
a basis of bideterminants, indexed by pairs of standard Young tableaux.
DeConcini-Eisenbud-Procesi [DEP] 
used these bideterminants to give an explicit filtration of  $K[X(i,j)]$,
as a $GL(n,K)$-bimodule,
whose factors are isomorphic to $L_\lambda\otimes {}_\lambda L$ where 
$L_\lambda$ and ${}_\lambda L$, are left and right Schur modules, respectively,
 corresponding
to partitions $\lambda$ into at most $n$ parts. At characteristic 0,
$L_\lambda$ and ${}_\lambda L$
are irreducible polynomial representations
of the general linear group $GL(n,K)$.

We wish to investigate the situation where $GL(n,K)$ is  replaced by
the orthogonal group $O(n,K)$. We assume that the characteristic of $K$ is
not 2. The bideterminants which are a basis of $K[X(i,j)]$ are not
linearly independent as functions on $O(n,K)$.
We show that the coordinate ring $K[O(n)]$ of $O(n,K)$
has a basis of
bideterminants $[S:T]$ where $S$ and $T$ are $O(n)$-standard in the
sense of King-Welsh [KW]. As far as we are aware, this is the first
known explicit basis of $K[O(n)]$.
We also find an explicit $O(n,K)$-bimodule filtration of $K[O(n)]$ 
with factors isomorphic to $L_O^\lambda\otimes {}^\lambda L_O$, where
$L_O^\lambda$ and ${}^\lambda L_O$, are analogues of 
left and right Schur modules, for the orthogonal group. 

Our first problem is to show that 
the $O(n)$-standard bideterminants generate $K[O(n)]$; to do this we
give a straightening algorithm. 
 We then give our filtration of
$K[O(n)]$, from which linear independence of the
$O(n)$-standard bideterminants can be deduced.

We first prove these results under the assumption that the field $K$
has characteristic 0. We then show that there are analogues of these
results where $K$ is replaced by the ring $\Z[1/2]$, and then by
base change, for fields of odd characteristic.

In [KW] a straightening method is given for an $O(n)$-module
(assuming $K=\C$) denoted $O^\lambda$ which is a factor module
of a submodule of $V^{\otimes m}$ where $V$ is the natural $O(n)$-module.
(This is defined in Section 6 below, where it is shown that $O^\lambda$
is isomorphic to $L^{\lambda}_O$.)
Straightening in [KW] is done using a modification of Berele's method [B]
(see also [Don])
for the symplectic group, using a suitable quotient of the tensor algebra.
Our straightening method is 
more complicated than that of [KW], as we need to use tableaux of different
shapes and sizes.
Linear independence of the basis in [KW] for
$O^\lambda$ is deduced from work
of Proctor [Pro].

We also give a basis of the space of homogeneous polynomial functions
of a given degree on the group of orthogonal similitudes. This will
be applied in future work to the orthogonal Schur algebra, as defined
by Doty [Dot].

In the case of the symplectic group, a basis involving bideterminants was
given by Oehms [O] for the space of homogeneous polynomial functions
of a given degree on
the group of symplectic similitudes. Our methods are quite different
from those of Oehms.

Filtrations for the coordinate ring of the symplectic group were also
considered by de Concini [dC]. Filtrations of this sort have been
studied by Donkin and by Kopppinen [K] for connected reductive
groups $G$. The orthogonal group is not connected, and our emphasis
is on explicit filtrations. For connected, simply connected semi-simple
groups $G$, the coordinate ring
$\C[G]$ has a canonical basis, due to Lusztig, [L], Chapter 29. This basis
is not
explicitly given.

The author would like to thank V.~Chernousov, A.~Pianzola, and the referee,
for comments especially concerning Section 7.

\section{Preliminaries}

Let $n\ge3$ be a positive integer
and let $m$ be the greatest integer $\le n/2$.
Define the ordered set
$\cI=\{\ov1<1<\ov2<2<\cdots<\ov m<m\}$ if $n$ is even, and
$\cI=\{\ov1<1<\ov2<2<\cdots<\ov m<m<0\}$ if $n$ is odd. 
Define
$$\ov{\bi}=i, \quad 1\le i\le m,\qquad \ov0=0.$$
Let $\{v_i:i\in\cI\}$ be the standard basis of $V=K^n$, and define
the symmetric bilinear form on $V$
\begin{equation}\label{bilform}
\left<\sum_{i\in\cI}x_iv_i,\sum_{i\in\cI}y_iv_i\right>
=\sum_{i\in\cI}x_iy_{\bi}=\sum_{i=1}^m(x_iy_{\bi}+x_{\bi}y_i)+x_0y_0
\end{equation}
where $x_0=0$ if $n$ is even.
The orthogonal group $O(n,K)$, which we will also denote by $O(n)$, is
defined as the subgroup of $GL(n,K)$ which preserves the 
form $\left<~,~\right>$  
given in (\ref{bilform}).

As usual, a {\em partition} $\lambda$ of a positive integer $r$
into $k$ parts is
given by writing $r$ as a sum $r=\lambda_1+\lambda_2+\cdots+\lambda_k$
of positive integers
where $\lambda_1\ge\lambda_2\ge\cdots\ge\lambda_k$.
We let $|\lambda|$ denote $r$, and we call $|\lambda|$ the {\em size} of
$\lambda$.
A {\em Young tableau} of shape $\lambda$ is a left justified array, having
$k$ rows; the $i$-th row has $\lambda_i$ entries from the set $\cI$.
We will denote the $i$-th column
of $T$ by $T_i$. 

The {\em conjugate} $\lambda'$ of a partition $\lambda$ is the partition
whose parts are the column lengths of a tableau of shape $\lambda$.

The {\em dominance} order on partitions of $r$ is defined by
$\lambda \unlhd \mu$ if $\lambda_1 + \cdots + \lambda_i \leq
\mu_1 + \cdots + \mu_i$ for all $i \in \{1,\ldots,k\}$,
where $\lambda$ has $k$ parts. (For this to make sense
we let $\mu_j=0$ if $j$ is greater than the number of parts of $\mu$.) For
example, if $\lambda=(2,2,1)$ and $\mu=(4,1)$, then $\lambda
\unlhd \mu$. 
We will use the following elementary result.
\begin{remark}\label{bottom}
If $T$ is a tableau of shape $\lambda$ and the tablea+u $S$ is obtained from
$T$ by moving an entry from the bottom of a column to the bottom
of a column to the left, then the shape $\mu$ of $S$ satisfies
$\mu\lhd\lambda$.
\end{remark}

A Young tableau 
 is called $GL(n)$-{\em standard} if it has at most $n$ rows, and
if the entries are non-decreasing from left to right across each row
and strictly increasing from top to bottom down each column. We will 
give the definition, due to [KW], of what is called an $O(n)$-{\em standard}  
Young tableau $T$. For $i=1,2,\ldots,m$ let $\alpha_i$ and 
$\beta_i$ be the number of entries less than or equal to $i$ in the first
and second columns, respectively, of $T$. Let $T(i,j)$ denote
the entry in row $i$ and column $j$ of $T$. We shall call $T$  
$O(n)$-{\em standard} if it is $GL(n)$-standard, 
if $\lambda'_1+\lambda'_2\le n$, 
and if in addition, for
each $i=1,2,\ldots,m$,
\begin{enumerate}\item[(OS 1)] $\alpha_i+\beta_i\le 2i$;
\item[(OS 2)] if $\alpha_i+\beta_i=2i$ with $\alpha_i>\beta_i$ and $T(\alpha_i,1)=i$
and $T(\beta_i,2)=\bi$ then $T(\alpha_i-1,1)=\bi$;
\item[(OS 3)] if $\alpha_i+\beta_i=2i$ with $\alpha_i=\beta_i$ ($=i$) and 
if $\bi,i$ occur in the $i$-th row of $T$, with $\bi$ in $T_1$
and $i$ in $T_b$ for some $b\ge 2$,
 then above the $i$ in $T_b$ there is an $\bi$.  
\end{enumerate}
In cases (2) and (3), the entry $i$ is said to be {\em protected} by
the existence of $\bi$ above it.

For each pair $(i,j)$ where $i$ and $j$ are in $\cI$, let $X(i,j)$
be an indeterminate, and let $X=(X(i,j))$ be the $n\times n$ matrix
whose rows and columns are indexed by $\cI$, and whose $(i,j)$-entry
is $X(i,j)$. Let $K[X(i,j)]$ be the polynomial ring in the $n^2$
indeterminates $X(i,j)$. 

The coordinate ring $K[O(n)]$ consists of the restrictions of the
functions in $K[X(i,j)]$ to the orthogonal group $O(n)$. (Strictly speaking,
a function in $K[O(n)]$ should be defined as a function on $O(n)$ given
by a polynomial in $K[X(i,j)]$
divided by some power of the determinant.
Since $\det=\pm1$ on $O(n)$, then $1/\det=\det$
as functions on $O(n)$; hence
all functions in $K[O(n)]$ are given by polynomials.)

Suppose that $S$ and $T$ are tableaux of the same shape $\lambda$.
Define
$[S_i:T_i]$ to be the determinant of the submatrix of $X$ whose rows
are indexed by the entries of the column $S_i$ and whose columns are indexed
by the entries of the column $T_i$. Define
the {\em bideterminant}, denoted $[S:T]$, to be the product
$$ [S:T]=\prod_{i=1}^{\lambda_1}\,[S_i:T_i].$$
We will sometimes use the notation
$$[S:T]_O$$
to denote the bideterminant $[S:T]$ considered as a function on $O(n)$.
We define the shape of $[S:T]$ or of $[S:T]_O$ to be $\lambda$, which is
the common shape of $S$ and $T$. We call $[S:T]_O$ an $O(n)$-standard
bideterminant if both $S$ and $T$ are $O(n)$-standard.

For an $n\times n$ matrix $A$, let $A^t$ denote the transpose of $A$.
\begin{remark}\label{transpose}The bideterminant $[S:T]$ evaluated at
$X^t$ is equal to the bideterminant $[T:S]$.
\end{remark}

For $1\le k\le n$, let
\begin{equation}\label{order}
j_1,j_2,\ldots,j_k\end{equation}
denote the first $k$ elements of the ordered set $\cI$.
For example, if $k=3$, then $j_1=\ov1$, $j_2=1$, $j_3=\ov2$.
The {\em basic} $\lambda$-tableau, denoted $T^\lambda$, is the tableau
having each entry in row $k$ equal to $j_k$.

We define the following partial order on the set of tableaux of shape 
$\lambda$. If $T\ne T'$, suppose that in the right-most column in which
there is a differing entry, the top-most entry $i$ of $T$ which differs
from an entry $i'$ of $T'$ in the same position satisfies $i<i'$; then
we say that $T\prec T'$.

By a {\em signed sum} of some quantities $x_1,x_2,\ldots, x_k$ we mean
a linear combination $\sum_{i=1}^k\epsilon_ix_i$ where each
$\epsilon_i\in\{-1,1\}$.

We say that $i,\bi$ {\em occur} in a 2-column tableau $T$ if $i$ is
an entry in the first column of $T$ and $\bi$ is an entry in the
second. To {\em delete} the pair $i,\bi$, remove $i$ from the first
column and $\bi$ from the second, and move up the other entries to
form a tableau.

We will let $\#A$ denote the size of a finite set $A$.

\section{$GL(n)$-Straightening}

Writing a bideterminant as a linear combination of $GL(n)$-standard ones is
often called {\em straightening}. We shall describe the straightening method
used by Mead [M].

Suppose that $S$ and $T$ are tableaux of the same shape $\lambda$ having
two columns, of lengths $k$ and $\ell$ respectively. Suppose that $S$ is 
column-increasing but not $GL(n)$-standard. Using Mead's notation for the
entries of the tableaux
(although Mead does not use tableaux or the term bideterminant) suppose that
$$[S:T]=\left[\begin{array}{cc}i_1&i'_1\\
                         \vdots&\vdots\\
                            i_\ell&i'_\ell\\
                         \vdots\\
                            i_k\end{array}:
           \begin{array}{cc}a_1&b_1\\
                         \vdots&\vdots\\
                            a_\ell&b_\ell\\
                         \vdots\\
                            a_k\end{array}\right]$$
Suppose that $i_j\le i'_j$ for $j<t$ but $i_t>i'_t$.

Let $H$ be the $(k+\ell)\times(k+\ell)$ matrix given by
$$H=\left(\begin{array}{cc}B_1&C\\D&B_2\end{array}\right)$$
where $B_1$ is the $k\times k$ submatrix of $X=(X(i,j))$ whose rows and columns
are indexed by $S_1$ and $T_1$ respectively,
$B_2$ is the
$\ell\times \ell$ submatrix of $X$
whose rows and columns are indexed by $S_2$ and $T_2$, respectively,
$C$ is the $k\times\ell$ matrix defined by
$$c_{p,q}=0\mbox{ if }p< t,\quad c_{p,q}=X(i_p,b_q)\mbox{ if }p\ge t$$
and $D$ is 
the $\ell\times k$ matrix given by
$$d_{p,q}=X(i'_p,a_q)\mbox{ if }p\le t,\quad d_{p,q}=0\mbox{ if }p> t.$$
Use Laplace expansion along multiple rows or columns 
(see, for example, [Pra, 2.4.1, p.~11]), in two ways, as follows:
expand $\det H$ by minors of the first $k$ columns, and
by minors of $\ell-1$ rows consisting of the first $t-1$ rows along with
the last $\ell-t$ rows. (An example will be given
below.)

The bideterminant $[S:T]$ occurs as one of the terms of the first
expansion; solve for $[S:T]$ equal to the negative of the signed sum
of the other bideterminants from the first expansion plus the signed
sum of bideterminants from the second expansion.

All the terms in the row-expansion involve bideterminants where the
tableaux have
2 columns, of length $k+1,\ell-1$, respectively, hence of shape
$\lhd\lambda$ by Remark \ref{bottom}.

In the column expansion, each term is a product of two minors, the
second of which has rows indexed by
$e_1,e_2,\ldots,e_{k},i'_{t+1},\ldots, i'_\ell$ where 
$\{e_1,e_2,\ldots,e_{k}\}$ vary over all possible subsets of
$\{i_t,\ldots,i_k,i'_1,\ldots,i'_t\}$.  The term $[U:T]$ with $U$ the
lowest (in the order $\prec$) is the one where $\{
e_1,e_2,\ldots,e_t\}=\{i'_1,\ldots,i'_t\}$; this gives the
bideterminant $[S:T]$.
All the
other bideterminants are of the form $[U:T]$ where $U\succ S$.

So we have the following result.

\begin{lemma}\label{mead}
Let $S$ and $T$ be two-column tableaux of the same shape $\lambda$,
where $S$ is column
increasing but not $GL(n)$-standard. Then (i)
$$[S:T]=\sum_Ua_U[U:T]+s$$ 
where each $U$ occurring in the sum is a tableau of shape $\lambda$
where $U\succ S$, each $a_U\in\{1,-1\}$, and $s$ is signed
sum of bideterminants  of shapes $\mu\lhd\lambda$.
(ii) The tableaux $U$ in the sum, and the signs $a_U$ depend only on $S$
and not on $T$. (iii) If $T=T^\lambda$, then $s=0$. 
\end{lemma}
\begin{proof}
We need only prove (ii) and (iii). For independence
of $T$, the row indices in the column expansion of $\det H$ depend only on $S$.
For (iii), if the rows of $T$ all have equal entries,
in the row expansion of $\det H$, 
each bideterminant we get
 has the form $[V:W]$ where the first column of
 $V$ has repeated entries and so
$[V:W]=0$.
\end{proof}
\begin{example}\label{mead-eg} Consider the bideterminant
$$[S:T]=\left[\begin{array}{cc}1&1\\
                               2&\ov2\\
                               3
                            \end{array}:
           \begin{array}{cc}\ov1&1\\
                             \ov2&2\\
                               3
                            \end{array}\right].$$
\end{example}
$$\mbox{Then }H=\left(
\bigskip\begin{array}{ccccc}
X(1,\ov1) & X(1,\ov2) & X(1,3) & 0 & 0\\
X(2,\ov1) & X(2,\ov2) & X(2,3) & X(2,1) & X(2,2) \\
X(3,\ov1) & X(3,\ov2) & X(3,3) & X(3,1) & X(3,2) \\
X(1,\ov1) & X(1,\ov2) & X(1,3) & X(1,1) & X(1,2) \\
X(\ov2,\ov1) & X(\ov2,\ov2) & X(\ov2,3) & X(\ov2,1) & X(\ov2,2) 
\end{array}\right).$$
For the expansion of $\det H$ along the first three columns, we get terms which
are a product of
two minors: in the first minor, we use the first row, we omit the fourth
row, since its first three entries are the same as those of the first row, 
and we have a choice of two
of rows 2, 3, and 5.
The second minor has the complementary rows.
 This gives
$$\det H=\left[\begin{array}{cc}1&1\\
                               2&\ov2\\
                               3\end{array}:
              \begin{array}{cc}\ov1&1\\
                               \ov2&2\\
                               3\end{array}\right]+
         \left[\begin{array}{cc}1&3\\
                                2&1\\
                               \ov2\end{array}:
              \begin{array}{cc}\ov1&1\\
                               \ov2&2\\
                               3\end{array}\right]-
         \left[\begin{array}{cc}1&2\\
                                3&1\\
                             \ov2\end{array}:
              \begin{array}{cc}\ov1&1\\
                               \ov2&2\\
                               3\end{array}\right].    
$$

The row expansion is along the first row, giving three terms:
$$\det H=        \left[\begin{array}{cc}1 & 1 \\	  
                                        2     \\	  
                                        3     \\	  
                                       \ov2   \end{array}:
                   \begin{array}{cc} \ov2 &\ov1 \\
                                        3       \\
                                        1       \\
                                        2       \end{array}\right]-
                 \left[\begin{array}{cc}1 & 1 \\
                                        2     \\
                                        3     \\
                                       \ov2   \end{array}:
              \begin{array}{cc}\ov1&\ov2\\
                               3\\
                               1\\
                               2\end{array}\right]+
                 \left[\begin{array}{cc}1 & 1 \\	  
                                       	2     \\	  
                                       	3     \\	  
                                        \ov2   \end{array}:
              \begin{array}{cc}\ov1&3\\
                               \ov2\\
                               1\\
                               2\end{array}\right].
$$
We switch rows to give standard tableaux. Then
$$\begin{aligned}
\left[\begin{matrix}1&1\\
                    2&\ov2\\
                    3\end{matrix}:
      \begin{matrix}\ov1&1\\
                    \ov2&2\\3
      \end{matrix}\right]&=
\left[\begin{matrix}1&1\\
                 \ov2&2\\
                 3\end{matrix}:
              \begin{array}{cc}\ov1&1\\
                               \ov2&2\\
                               3\end{array}\right]
-\left[\begin{array}{cc}1&1\\
                      \ov2&3\\
                      2\end{array}:
              \begin{array}{cc}\ov1&1\\
                               \ov2&2\\
                               3\end{array}\right]
\\&+   
\left[\begin{array}{cc}1&1\\
                               \ov2\\
                               2\\
                               3\end{array}:
              \begin{array}{cc}1&\ov1\\
                               2\\
                               \ov2\\
                               3\end{array}\right]-
        \left[\begin{array}{cc}1&1\\
                               \ov2\\
                               2\\
                               3\end{array}:
              \begin{array}{cc}\ov1&\ov2\\
                               1\\
                               2\\
                               3\end{array}\right]+
        \left[\begin{array}{cc}1&1\\
                               \ov2\\
                               2\\
                               3\end{array}:
              \begin{array}{cc}\ov1&3\\
                               1\\
                               \ov2\\
                               2\end{array}\right].
\end{aligned}$$

We have the following result on $GL(n)$-straightening.
\begin{theorem}\label{GL(n)-straight}Suppose that $S$ and $T$
are tableaux of the same shape $\lambda$. Then
$$[S:T]=\sum_Ub_U[U:T]+s$$
where the tableaux $U$ in the sum are $GL(n)$-standard, each $b_U\in K$
and is independent of $T$, and $s$ is a signed sum of bideterminants
of shapes $\mu\lhd\lambda$. 
\end{theorem}
\begin{proof}This is proved, using Lemma \ref{mead},
 by induction on the order $\lhd$
and downward induction on $\prec$.  
\end{proof}
Using Remark \ref{transpose}
there is an analogous version where the roles of $S$ and $T$ are
interchanged, writing $[T:S]$ as a linear combination of $[T:U]$
plus bideterminants of shape smaller than $\lambda$ in the $\lhd$ order. 
This gives the spanning part of the following
well-known result.
For linear independence, which we do not need, see [M], [DRS] or [DEP].

\begin{theorem}\label{GLn-straight-basis}
The bideterminants $[S:T]$ where both $S$ and $T$ are
$GL(n)$-standard are a basis of the ring $K[X(i,j)],1\le i,j\le n$.
\end{theorem}

In Sections 7 and 8 we will replace $K$ with the ring $\Z$ or $\Z[1/2]$.
In the statement of Theorem \ref{GL(n)-straight}, the 
coefficients $b_U$  come from the coefficients $a_U$ of Lemma \ref{mead},
and these are $\pm1$. So we have the following.

\begin{theorem}\label{GLn-straight, any R}
The bideterminants $[S:T]$ where both $S$ and $T$ are
$GL(n)$-standard are a basis of the ring $R[X(i,j)],1\le i,j\le n$,
where $R$ is any commutative ring.
\end{theorem}
As with Theorem \ref{GLn-straight-basis} we will not need to
use linear independence.
We will also use the following.
\begin{lemma}\label{one-switch}
Suppose that $S$ and $T$ are tableaux of the same shape $\lambda$ having
two columns, where $S$ is $GL(n)$-standard.
 Suppose that
$\bi,i$ occur in $S$, both in the $t$-th row,
and let $S^*$ be the result of replacing this pair
$\bi,i$ with $i,\bi$. 
Then 
$$[S^*:T]=[S:T]+\sum_Ua_U[U:T]+s$$
where each $U$ is a $\lambda$-tabeau such that $S^*\prec U$,
each $a_U$ is $\pm1$, and $s$ is signed sum of bideterminants
of shape $\mu$ where $\mu\lhd\lambda$. The $U$ and $a_U$ occurring in
the equation are independent of $T$.
\end{lemma}
\begin{proof}
In Mead's method,
in the column expansion, the lowest $U$ is $S$, and the second lowest
is $S^*$. \end{proof}
This is illustrated in Example \ref{mead-eg} above,
where $[S:T]$ is written as a signed sum of standard bideterminants,
the first of which is $[S^*:T]$ where $S^*$ is obtained from 
$S$ by replacing 
$2,\ov2$ in the second row by $\ov2,2$, and
and the second is $[U:T]$ where $U\succ S^*$.

\section{Main Technical Result}
As functions on $O(n)$, we have
$$\sum_{i\in\mathcal I}X(i,j)X(\bi,k)=\delta_{j\ov k},\quad j,k\in\cI.$$
In bideterminant notation, 
$$X(i,j)X(\bi,k)=[i\ \bi: j\ k ] .$$
So we have, on $O(n)$,
$$\sum_{i\in\mathcal I}\,
[i\ \bi : j\ k]=\delta_{j\ov k}.$$
We will need more general versions of this. For example, suppose that
$S$ is a tableau having 2 columns, of lengths $e$ and $f$,
and that $T$ is a tableau of 2 columns, of lengths $e+1$ and $f+1$.
It will follow from our results below that if 
$i,\bi$ do not occur in $T$ for any $i\in\cI$, then
$$\sum_{i\in\cI}\,
\left[\begin{array}{cc}
  i& \bi\\ \multicolumn{2}{c}S
\end{array}
: T\right]
=0$$
on $O(n)$.
More generally, if $i,\bi$ do occur in $T$ for some $i\in\cI$, then
$$\sum_{i\in\cI}\,
\left[\begin{array}{cc}
  i& \bi\\ \multicolumn{2}{c}S
\end{array}
: T\right]
=
\sum_{T'}\pm\,[S:T']$$
on $O(n)$, 
where the sum is over all tableaux $T'$ obtained by
deleting a pair $i,\bi$ from $T$ if $i,\bi$  occur in $T$.

For a more general version of this, we will have not just one
variable $i$ in the summation, but $a$ variables $i_1,i_2,\ldots,
i_a$, each of which will vary over $\cI$, except that each $i_j$ will
not be allowed to assume one of $c$ fixed values, where $0\le c<a$.

Suppose that $S$ and $T$ are tableaux of the same shape, having two columns.
Suppose that $a$ is a positive integer, less than or equal to the length
of the second column of $S$ (and of $T$). Suppose that $C$ is a 
(possibly empty) subset
of $\cI$ having $c$ elements, where $c<a$. 
Let $S_0$ be the tableau obtained from $S$ by deleting its first $a$ rows.
We want to calculate the sum
\begin{equation}\label{sum2}L=\sum_{i_1,\ldots,i_a\in\cI-C}
\left[
\begin{array}{cc}i_1&\bi_1\\
              i_2&\bi_2\\
               \vdots&\vdots\\
              i_a&\bi_a\\
              \multicolumn{2}{c}{S_0}
               \end{array}:T\right].
\end{equation}
We allow repeated indices in the sum, although a bideterminant
having repeated indices $i_s$, is 0.

Suppose that $d$ is an integer, $1\le d\le a$, and suppose
that $E=\{i_1,\ldots,i_d\}$ is a subset of $\cI$ such that
$i,\bi$ occur in $T$ for all $i\in E$. Let
$$(T,E)$$
be the tableau obtained by deleting the pairs $i,\bi$
from $T$ for each $i\in E$.
Let $\cE_d$ denote the set of all subsets $E$ of $\cI$ of size
$d$ such that $\ov\imath,i$ occur in $T$ for all $i\in E$.

We next define $\cS_d$ by the equation 
\begin{equation}\label{Sd}
\cS_d=\sum_{E\in\cE_d}\,
\sum_{\substack{i_1<\cdots< i_{a-d}\\ i_j\in C}}
(-1)^{a-d}
\left[
\begin{array}{cc}
i_1&\bi_1\\
              i_2&\bi_2\\
               \vdots&\vdots\\
              i_{a-d}&\bi_{a-d}\\
              \multicolumn{2}{c}{S_0}
               \end{array}:(T,E)\right].
\end{equation}
If for a given $d$, $C$ has fewer than $a-d$ elements, then the second sum
is the empty sum, and $\cS_d=0$. In particular, if $C$ is empty,
then each $\cS_d=0$ except for $d=a$.
\bigskip

\begin{lemma}\label{main}
With the above notation, the sum $L$ in (\ref{sum2}),  
as function on $O(n)$, is given by
$$L=a!\sum_{d=1}^a\cS_d.$$
\end{lemma}

Before giving the proof, we have an example.
Consider
$$
\sum_{\substack{i_1,i_2,i_3\in\cI\\i_1,i_2,i_3\ne\ov4,5}}
\left[
\begin{array}{cc}
i_1&\bi_1\\
i_2&\bi_2\\
i_3&\bi_3\\
           7&8\\
           9
\end{array}:
\begin{array}{cc}%
\ov1&1\\
\ov2&2\\
\ov3&3\\
4&5\\6&\\
\end{array}
\right].$$

Then $a=3$, $C=\{\ov4, 5\}$, $c=\#C=2$. First suppose that $d=1$.
Since $a-d=2$, there is a unique subset of $C$ of size 2, so we replace
$(i_1,\bi_1),(i_2,\bi_2),(i_3,\bi_3)$ with the two pairs $(\ov4,4),(5,\ov5)$. 
There are three subsets $E_1$.
$$\cS_1
=
\left[\begin{array}{cc}
\ov4&4\\
5&\ov5\\
 7&8\\
 9
\end{array}:
\begin{array}{cc}
\ov2 & 2\\
\ov3 & 3\\
 4&5\\6&\\
\end{array}
\right]
+\left[\begin{array}{cc}
\ov4&4\\
5&\ov5\\
 7&8\\
 9
\end{array}:
\begin{array}{cc}
\ov1 & 1\\
\ov3 & 3\\
 4&5\\6&\\
\end{array}
\right]
+\left[\begin{array}{cc}
\ov4&4\\
5&\ov5\\
 7&8\\
 9
\end{array}:
\begin{array}{cc}
\ov1 & 1\\
\ov2 & 2\\
 4&5\\6&\\
\end{array}
\right].
$$
When $d=2$, we  replace the pairs $(i_1,\bi_1)$, $(i_2,\bi_2)$,
$(i_3,\bi_3)$ from $S$ with $a-d=1$ pair $i,\bi$ where $i\in C$. There
are two ways to do this, since $C$ has two elements. There are three subsets
$E_2$. So
$$\cS_2=-\sum_{i=\ov4,5}\sum_{j=1}^3
\left[\begin{array}{cc}
i & \bi\\
7&8\\9\end{array}:
\begin{array}{cc}
\bj & j\\
4 &5\\6\end{array}
\right]$$
When $d=3$, there is only one $E_3$, namely $\{1,2,3\}$. Since
$a-d=0$,  
$$\cS_3=\left[\begin{array}{cc}
7&8\\9\end{array}:
\begin{array}{cc}
4 &5\\6\end{array}
\right]$$

\begin{proof}
Suppose that the entries of the two columns of $T$ are
$$g_1,g_2,\ldots,g_e,\qquad
h_1,h_2,\ldots,h_f,\mbox{ respectively.}$$
Suppose further that the entries of the two columns of $S_0$ are
$$j_{a+1},\ldots,j_e,\qquad k_{a+1},\ldots,k_f\mbox{~~respectively.}$$
Let $P$ be the permutations of the symbols $g_1,g_2,\ldots,g_e$,
and $Q$ the permutations of the symbols $h_1,\ldots,h_f$.

Then the sum $L$ in (\ref{sum2}) is given by
\begin{multline*}
L=\sum_{i_1,\ldots,i_a\in\cI-C}
\sum_{\sigma\in P}\sgn(\sigma) 
\prod_{s=1}^{a}X\bigl(i_s,\sigma(g_s)\bigr)
\prod_{s=a+1}^{e}X\bigl(j_s,\sigma(g_s)\bigr)\\
\times\sum_{\tau\in Q}\sgn(\tau)
\prod_{t=1}^{a}X\bigl(\bi_t,\tau(h_t)\bigr)
\prod_{t=a+1}^{f}X\bigl(k_t,\tau(h_t)\bigr).
\end{multline*}
Define 
$$Y_{\sigma,\tau}=\prod_{s=a+1}^{e}X\bigl(j_s,\sigma(g_s)\bigr)
\prod_{t=a+1}^{f}X\bigl(k_t,\tau(h_t)\bigr).$$
Then we have
$$L=\sum_{\sigma\in P} \sum_{\tau\in Q}\sgn(\sigma)\sgn(\tau)
\prod_{s=1}^{a}
\left(\sum_{i_s\in \cI-C}X\bigl(i_s,\sigma(g_s)\bigr)
X\bigl(\bi_s,\tau(h_s)\bigr)\right)Y_{\sigma,\tau}.$$
Since $\sum_{i\in\cI}X(i,j)X(\bi,k)=\delta_{j,\ov k}$ on $O(n)$, then
as functions on $O(n)$, for each $s$ in $\{1,\ldots,a\}$,
\begin{equation}\label{main1}\sum_{i_s\in\cI-C}
X\bigl(i_s,\sigma(g_s)\bigr)X\bigl(\bi_s,\tau(h_s)\bigr)
=\delta_{\sigma(g_s),\ov{\tau(h_s)}}
-\sum_{i_s\in C}X\bigl(i_s,\sigma(g_s)\bigr)
X\bigl(\bi_s,\tau(h_s)\bigr).
\end{equation}
Hence as a function on $O(n)$, $L$ is equal to
$$L'=\sum_{\sigma\in P} \sum_{\tau\in Q}\sgn(\sigma)\sgn(\tau)
\prod_{s=1}^{a}\left(\delta_{\sigma(g_s),\ov{\tau(h_s)}}
-\sum_{i_s\in C}X\bigl(i_s,\sigma(g_s)\bigr)
X\bigl(\bi_s,\tau(h_s)\bigr)
\right)Y_{\sigma,\tau}.$$
To simplify the notation somewhat, fix $(\sigma,\tau)\in P\times Q$, and let
$$\delta(s)=\delta_{\sigma(g_s),\ov{\tau(h_s)}},\quad
Z(s,i_s)=X\bigl(i_s,\sigma(g_s)\bigr)X\bigl(\bi_s,\tau(h_s)\bigr).$$
With this notation we have, for fixed $\sigma,\tau$,
$$\prod_{s=1}^{a}\left(\delta_{\sigma(g_s),\ov{\tau(h_s)}}
-\sum_{i_s\in C}X\bigl(i_s,\sigma(g_s)\bigr)
X\bigl(\bi_s,\tau(h_s) \bigr)\right)=
\prod_{s=1}^a\left(\delta(s)-\sum_{i_s\in C}Z(s,i_s)\right).$$
Expand the product
$$\prod_{s=1}^a\left(\delta(s)-\sum_{i_s\in C}Z(s,i_s)\right).$$
To do this, for each $s$ select a term from the factor
$\left(\delta(s)-\sum_{i_s\in C}Z(s,i_s)\right)$, namely 
$\delta(s)$ or $-Z(s,i_s)$ for some $i_s\in C$; multiply all the chosen
terms, and sum over all possible choices of terms.
For a given choice,
let $D$ be the set of all $s$ for which we pick $\delta(s)\ne0$. If
$s\in D$, since
$$0\ne\delta(s)=\delta_{\sigma(g_s),\ov{\tau(h_s)}}$$ 
then
$$\sigma(g_s)=\ov{\tau(h_s)}.$$
Let $\cD$ be the set of all subsets $D$, including the empty subset, 
of $\{1,\ldots,a\}$ such that
$\delta(s)\ne0$ for all $s\in D$. For a given $D\in\cD$,
let 
$$\ell=\ell(D)=a-\#D$$
 and let $D'$ be the complement of $D$ in
$\{1,2,\ldots,a\}$. Then 
\begin{equation}\label{e-main2}
\prod_{s=1}^a\left(\delta(s)-\sum_{i_s\in C}Z(s,i_s)\right)
=\sum_{D\in\cD}(-1)^{\ell(D)}\sum_{(i_1,\ldots,i_\ell)\in C^\ell}
\prod_{s\in D'}Z(s,i_s).
\end{equation}
Write $\cD=\cD(\sigma,\tau)$.
Then as functions on $O(n)$, $L$ is equal to 
$$L''=\sum_{\substack{\sigma\in P\\ \tau\in Q}}
\sgn(\sigma)\sgn(\tau)
\sum_{D\in\cD(\sigma,\tau)}(-1)^{\ell(D)}\sum_{(i_1,\ldots,i_\ell)\in C^\ell}\,
\prod_{s\in D'}X(i_s,\sigma(g_s))X(\bi_s,\tau(h_s))
Y_{\sigma,\tau}.$$

We want to interchange the order of the first two summations.
Let $\cP$ denote the set of all subsets of $\{1,\ldots,a\}$.
For a given subset $D$ of $\cP$, define
$$G(D)=\{(\sigma,\tau)\in P\times Q:
\sigma(g_s)=\ov{\tau(h_s)}\mbox{ for all }s\in D\}.$$
If $D$ is empty, then $G(D)=P\times Q$.
Then
$$L''=\sum_{D\in\cP}(-1)^{\ell(D)}\!\!\!\!\sum_{(\sigma,\tau)\in G(D)}
\sgn(\sigma)\sgn(\tau)
\!\!\!\sum_{(i_1,\ldots,i_\ell)\in C^\ell}\,
\prod_{s\in D'}X(i_s,\sigma(g_s))X\bigl(\bi_s,\tau(h_s)\bigr)
Y_{\sigma,\tau}.$$
For a given subset $D$ of $\cP$, let
$$L''(D)=\sum_{(\sigma,\tau)\in G(D)}\sgn(\sigma)\sgn(\tau)
\sum_{(i_1,\ldots,i_\ell)\in C^\ell}\,
\prod_{s\in D'}X(i_s,\sigma(g_s))X\bigl(\bi_s,\tau(h_s)\bigr)
Y_{\sigma,\tau}.$$
Hence
$$L''=\sum_{D\in\cP}(-1)^{\ell(D)} L''(D).$$
Let $\cP_d$ be the set of subsets of $\{1,\ldots,a\}$ of size $d$, and
let 
$$L''_d=\sum_{D\in\cP_d}(-1)^{a-d} L''(D);\quad\mbox{then}\quad
L''=\sum_{d=0}^a L''_d.$$
For a fixed $D\in\cP_d$, let $\bolde$ be a $d$-tuple
$\{e_1,e_2,\ldots,e_d\}$ of distinct elements of $\cI$. Let
$$P(D,\bolde)=\{\sigma\in P:
\sigma(g_s)=e_s,s\in D\},~
Q(D,\bolde)=\{\tau\in Q:\tau(h_s)=\ov e_s,s\in D\}.$$
If $D$ is empty, then $\bolde$ is a 0-tuple, and then 
$P(D,\bolde)=P$, and $Q(D,\bolde)=Q$.
Let
$$L''(D,\bolde)=\sum_{\substack{\sigma\in P(D,\bolde)\\ \tau\in Q(D,\bolde)}}
\sgn(\sigma)\sgn(\tau)
\sum_{(i_1,\ldots,i_\ell)\in C^\ell}\,
\prod_{s\in D'}X(i_s,\sigma(g_s))X(\bi_s,\tau(h_s))
Y_{\sigma,\tau}.$$
Let $\cF_d$ denote the set of all $d$-tuples of distinct elements of
$\cI$;
then
$$L''(D)=\sum_{\bolde\in\cF_d}L''(D,\bolde).$$
For a given $D$ and $\bolde$, consider a fixed $\ell$-tuple ($\ell=\ell(D)$) 
$$\boldi=(i_1,\ldots,i_\ell)\in C^{\ell}.$$
For this $\boldi$ let
$$L''(D,\bolde,\boldi)=
\sum_{(\sigma,\tau)\in P(D,\bolde)\times Q(D,\bolde)}\sgn(\sigma)\sgn(\tau)
\prod_{s\in D'}X\big(i_s,\sigma(g_s)\bigr)X\bigl(\bi_s,\tau(h_s)\bigr)
Y_{\sigma,\tau},$$
so 
$$L''(D,\bolde)=\sum_{\boldi\in C^\ell}L''(D,\bolde,\boldi).$$
Since 
$$Y_{\sigma,\tau}=\prod_{s=a+1}^{e}X\bigl(j_s,\sigma(g_s)\bigr)
\prod_{t=a+1}^{f}X\bigl(k_t,\tau(h_t)\bigr)$$
then 
$L''(D,\bolde,\boldi)$ can be factored as the product 
$M(D,\bolde,\boldi)N(D,\bolde,\boldi)$ where
$$M(D,\bolde,\boldi)=\sum_{\sigma\in P(D,\bolde)}\sgn(\sigma)
\prod_{s\in D'}X(i_s,\sigma(g_s))
\prod_{s=a+1}^{e}X(j_s,\sigma(g_s)),$$
$$N(D,\bolde,\boldi)=\sum_{\tau\in Q(D,\bolde)}\sgn(\tau)
\prod_{s\in D'}X\bigl(\bi_s,\tau(h_s)\bigr)
\prod_{s=a+1}^{f}X\bigl(k_s,\tau(h_s)\bigr).$$

Up to a sign, $M(D,\bolde,\boldi)$ is the determinant of the
submatrix of $X$ whose rows are indexed by $i_1,\ldots, i_\ell$ 
(where $\boldi=(i_1,\ldots,i_\ell)$) and $j_{a+1},\ldots,j_e$,
and whose columns are indexed by
the first column of $T$
with  the set $E$ of entries of $\bolde$ deleted. Further, up to the same sign,
$N(D,\bolde,\boldi)$ is the determinant of the
matrix whose rows are indexed by 
$\bi_1,\ldots \bi_\ell$ and
$k_{a+1},\ldots,k_f$, and whose columns are indexed by
the second column of $T$,
with  $\ov E$ deleted. 
The product of the two signs is 1, so
\begin{equation}\label{L'' }
L''(D,\bolde,\boldi)= 
\left[
\begin{array}{cc}i_1&\bi_1\\                 i_2&\bi_2\\
               \vdots&\vdots\\
               i_{\ell}&\bi_{\ell}\\
              \multicolumn{2}{c}{S_0}
               \end{array}:(T,E)\right].
\end{equation}
This is one of the terms in equation (\ref{Sd}).
Note that the right side of (\ref{L'' })
does not depend on $D$. Since there are $\binom ad$
subsets of $\{1,\ldots,a\}$ of size $d$, then for any $d$, pick
some $D_d\in\cP_d$, and then we have
$$L''_d=\binom ad\sum_{\bolde\in\cF_d} (-1)^{a-d}L''(D_d,\bolde).$$

Also, the right side of equation (\ref{L'' })
does not depend on the ordering of the $d$-tuple $\bolde$,
just on the set $E$ of its entries. Write $L''(D,E)=L''(D,\bolde)$.
Each subset $E$ of $\cI$ of size $d$ gives us $d!$ distinct $d$-tuples in $\cF_d$
by permuting the elements. Recall that $\cE_d$ is the set of all subsets
of $\cI$ of size $d$.
So we have
$$L''_d=\binom add!\sum_{E\in\cE_d} (-1)^{a-d}L''(D_d,E).$$

If $\boldi$ has two equal entries, then the right side of
(\ref{L'' }) is 0, so $L''(D,E,\boldi)=0$;
hence we assume that all the entries of $\boldi$ are distinct.
We did not exclude the possiblity that $D$ is the empty set. If it is,
then $\ell(D)=a$, so $i_1,\ldots,i_\ell$ are $a$ elements of $C$,
but we are assuming that $\#C<a$. So if $D$ is the empty set, then $\boldi$
has repeated entries, which we have excluded. Hence
 we may assume that $d\ge 1$,
and it follows that
$$L''=\sum_{d=1}^aL''_d.$$

If we permute the elements of $\boldi$, giving say $\boldi'$, then
$M(D,\bolde,\boldi')$ and $N(D,\bolde,\boldi')$ are, up to the same sign,
equal to $M(D,\bolde,\boldi)$ and $N(D,\bolde,\boldi)$, respectively.
So $L''(D,E,\boldi')=L''(D,E,\boldi)$. Let
$$C(\ell)=\{(i_1,\ldots,i_\ell):i_1<\ldots<i_\ell, i_j\in C, j=1,\ldots,\ell\}.
$$
Then 
$$L''(D,E)=\sum_{\boldi\in C^{\ell}}L''(D,E,\boldi)
=\ell!\sum_{\boldi\in C(\ell)}L''(D,E,\boldi).$$
Since $\ell=a-d$ then
$$L''_d=\binom add!(a-d)!\sum_{E\in\cE_d}
\sum_{\boldi\in C(\ell)}L''(D,E,\boldi).$$

Since $\binom ad d!(a-d)!=a!$, then
$$L''_d=\sum_{E\in\cD_d}L(D,E)=a!
\sum_{d=1}^a\cS_d.$$
Since $L=L''$ as functions on $O(n)$, the proof is complete.
\end{proof}

In the definition of $L$ in (\ref{sum2}) 
we allowed all possible $i_1,\ldots,i_a\in\cI-C$. We 
define 
\begin{equation}\label{sum3}
\tilde L=
\sum_{\substack{i_1,\ldots,i_a\in\cI-C\\
                  i_1<\cdots<i_a}}
\left[
\begin{array}{cc}i_1&\bi_1\\
              i_2&\bi_2\\
               \vdots&\vdots\\
              i_a&\bi_a\\
              \multicolumn{2}{c}{S_0}
               \end{array}:T\right].
\end{equation}

\begin{lemma}\label{main2}
If $K$ has characteristic 0, then as functions on $O(n)$,
$\tilde L$ is equal to $\sum_{d=1}^a\cS_d$.
\end{lemma}
\begin{proof}
As we have seen, repeated indices contribute nothing to the sum $L$.
The result now follows from Lemma \ref{main} by dividing by $a!$.
\end{proof}
It will follow from our results in Section 7 that this result also 
holds at positive characteristic. However,
because of our characteristic zero assumption in this last lemma,
we will need to assume that $K$ has charactersistic 0 for most of our
results until the end of Section 6.

In Lemmas \ref{main} and \ref{main2}, we sum over $E\in \cE_d$; if 
for all $i\in \cI$, $i,\bi$ do not occur in $T$, then
each $\cE_d$ is empty, and $\cS_d=0$.
We then have the following.
\begin{lemma}\label{empty}
Suppose that for all $i\in\cI$, $i,\bi$ do not occur in $T$.
Then as functions on $O(n)$, 
$$\sum_{\substack{i_1,\ldots,i_a\in\cI-C\\
                  i_1<\cdots<i_a}}
\left[\begin{array}{cc}i_1&\bi_1\\
              i_2&\bi_2\\
               \vdots&\vdots\\
              i_a&\bi_a\\
\multicolumn{2}{c}{S_0}
               \end{array}:T\right]
=0.$$
\end{lemma}

A first application is that on $O(n)$, any bideterminant
is, up to sign, equal to a bideterminant
of shape $\lambda$, where $\lambda_1'+\lambda_2'\le n$.
We will prove this now in the two-column case.

\begin{lemma}\label{one-col} Supppose that $[S:T]$ is a bideterminant
whose shape consists of a single column. Let 
$\ov S'$ denote the column-increasing 
tableau whose entries are $\{\bi\in\cI:i\notin S\}$.
Then as functions on $O(n)$, 
\begin{equation}\label{e-one-col}
[S:T]=\pm\det\cdot[\ov S':\ov T'].
\end{equation}
\end{lemma}
\begin{proof} We assume that each of $S$ and $T$ has
no repeated entries. Suppose that $S$ has $k$ entries, and let $a=n-k$.
Let $V$ be the one-column tableau of length $n$ whose entries are
all the elements of $\cI$.
Consider the  sum
\begin{equation}\label{sum4}
s=\sum_{\substack{i_1,\ldots,i_a\in\cI\\i_1<\cdots<i_a}}
\left[\begin{array}{cc}i_1&\bi_1\\
              i_2&\bi_2\\
               \vdots&\vdots\\
               i_{a}&\bi_a\\
                      S \end{array}: V~\ov T'\right].
\end{equation}
Apply Lemma \ref{main2}. Since the set $C$ is empty, each $\cS_d$ is
empty except for $d=a$. We delete $a$ pairs $i,\bi$ from the tableau 
with columns $V$, $\ov T'$; this means we delete all the entries of $\ov T'$,
and what remains of $V$ is $T$. So we get, as functions on $O(n)$,
$$s=\pm[S: T].$$
On the other hand, unless $i_1,\ldots,i_a$  
all do not occur in $S$,
$$ \left[\begin{array}{c}i_1\\
              i_2\\
               \vdots\\
              i_{a}\\
                      S \end{array}: V     \right]=0.$$
So all terms in the sum (\ref{sum4}) are 0 except when
 $i_1,i_2,\ldots,i_a$ are the complementary entries to $S$, in which case
$$\left[\begin{array}{c}i_1\\
              i_2\\
               \vdots\\
              i_{a}\\
                      S \end{array}: V     \right]=\pm\det,\quad
\left[\begin{array}{cc}i_1&\bi_1\\
              i_2&\bi_2\\
               \vdots&\vdots\\
               i_{a}&\bi_a\\
                      S \end{array}: V~\ov T'\right]=
\pm\det\cdot[\ov S':\ov T'].
 $$
So we get $[S:T]=\pm\det\cdot[\ov S':\ov T']$, on $O(n)$, as desired.
\end{proof}
\begin{lemma}\label{tilde-lambda}
Suppose that $S$ and $T$ each have shape $\lambda$, where each
of $S$ and $T$ have two columns. Suppose that $\lambda_1'+\lambda_2'>n$.
Then $[S:T]$ is equal, up to sign, to a bideterminant $[\tilde S:\tilde T]$
where $\tilde S$ and $\tilde T$ have shape $\tilde\lambda$, where
each of $\tilde S,\tilde T$ have two columns, and 
$\tilde\lambda_1'+\tilde\lambda_2'<n$.
\end{lemma}
\begin{proof} Apply Lemma \ref{one-col}
to each of $[S_1:T_1]$ and $[S_2:T_2]$, giving, on $O(n)$
\begin{equation}\label{eqn-tilde-lambda}
[S:T]=\pm\det{}^2\cdot[\bar S_1':\bar T_1'][\bar S_2':\bar T_2'].
\end{equation}
If $g\in O(n)$, then $\det(g)=\pm1$, 
so we may delete the factor $\det^2$.
Let $[\tilde S:\tilde T]$ be the bideterminant where
 $\tilde S_1=\bar S'_2, \tilde S_2=\bar S'_1,
 \tilde T_1=\bar T'_2, \tilde T_2=\bar T'_1$. 
The shape of $[\tilde S:\tilde T]$
is $\tilde\lambda$ where $\tilde\lambda'_1=n-\lambda'_2$,
$\tilde\lambda'_2=n-\lambda'_1$. Then $[S:T]=\pm[\tilde S:\tilde T]$.
Since $\lambda_1'+\lambda_2'> n$ then
$$\tilde\lambda'_1+\tilde\lambda'_2=2n-\lambda_1'-\lambda_2'< n.$$
This completes the proof.
\end{proof} 

\section{$O(n)$-straightening}
We assume that $K$ has characteristic 0.
The following three lemmas are replacements of Lemma 3.6, 3.7, and 3.8
of [KW]. Our situation is more complicated than in [KW], as we need tableaux
of different shapes and sizes. 

We first deal with violations of condition (OS 1) in the definition of
$O(n)$-standard in Section 2. 
\begin{lemma}\label{3.6}
Suppose that $S$ and $T$ are $\lambda$-tableaux having two columns, and
that $S$ is not $O(n)$-standard, in that $\alpha_j+\beta_j>2j$ for some $j$.
Then 
\begin{equation}\label{equn3.6}[S:T]=\sum_{U\in\cS,U\ne S} -[U:T] + s\end{equation}
where $\cS$ is a set of tableaux 
all of shape $\lambda$ such that 
 $S\prec U$, and $s$ is 
a signed sum of bideterminants of 
shapes $\mu$ with $|\mu|<|\lambda|$. The set $\cS$ does not depend on $T$.
If $T=T^\lambda$, then $s=0$.
\end{lemma}
\begin{proof}
Let $\cJ=\{i\in\cI:i\le j\}$, and
$\cK=\{i\in\cI:i> j\}$; then $\#\cJ=2j$.
For a column $S_k$ of $S$, write $i\in S_k$ if $i$ occurs as an entry
in column $S_k$.
Let 
\begin{eqnarray*}
A&=&\{i\in\cJ:i\in S_1\mbox{ and }\bi\in S_2\}\\ 
B&=&\{i\in\cJ:i\in S_1\mbox{ or }\bi\in S_2\mbox{ but not both}\}\\
C&=&\{i\in\cJ:i\notin S_1, \bi\notin S_2\}.
\end{eqnarray*}
Then $\cJ$ is the disjoint union of $A$, $B$, and $C$.
Let $a=\#A,b=\#B,c=\#C$. Then $a+b+c=\#\cJ=2j$, and 
$\alpha_j+\beta_j=2a+b$. Since we are assuming that
$\alpha_j+\beta_j>2j$, then $2a+b>2j$, so $a>c$.

 Let $\cS$ be the set of all tableaux obtained from $S$ as follows.
For each $i\in A$, replace the pair $i,\bi$ in $S$ by a pair $i',\bi'$
where $i' \in \cI-C$; do this in all possible ways, where the
replacements are strictly increasing down the first column of
the tableau.

For each $U\in\cS$, let $U'$ be the tableau obtained 
from $U$ by rearranging
the elements in each column so that the replacements $i,i'$ occur in
the first $a$ rows. Then 
$$\sum_{U\in\cS}[U:T]=\pm \sum_{U\in\cS}[U':T]$$
and this last sum has the form (\ref{sum3}). 
From Lemma \ref{main2}, the sum $\sum_{U\in\cS}[U:T]$
is equal to a signed sum $s$ of bideterminants having shapes $\mu$ such that
$|\mu|<|\lambda|$, and from Lemma \ref{empty}, $s=0$ if $T=T^\lambda$.
If for some $U\in\cS$ a replacement $i'$ is in $B$, then either
$i'$ occurs twice in $U_1$ or $\bi'$ occurs twice in $U_2$, so
the bideterminant $[U:T]=0$. 
So the only non-zero bideterminants in the sum $\sum_{U\in\cS}[U:T]$
are those for which all the replacements come from $A\cup\cK$,
since the replacements cannot come from $C$. One of these
is the original $[S:T]$, where each pair $i,\ov\imath$ is replaced by itself.
Any other $U\in\cS$ has at least one replacement $i'$ coming from 
$\cK$, so $i'>j$
and therefore $S\prec U$. Then
$$[S:T]=-\sum_{U\in\cS,U\ne S}[U:T]+s.$$
and for each $U$ where $[U:T]$ occurs in this sum, $S\prec U$.
The tableaux $U$ in this sum depend only on $S$ and not on $T$.
\end{proof}

\noindent Example. Consider the bideterminant as a function on $O(6)$
$$[S:T]=\left[\begin{array}{cc}\ov1&\ov2\\ \ov2&2\\2\end{array}:
              \begin{array}{cc}1&2\\\ov2&3\\\ov3\end{array}\right].$$
For the tableau $S$, $\alpha_2=3, \beta_2=2, \alpha_2+\beta_2>4$.
For $j=2$, $A=\{\ov2,2\},B=\{\ov1\},C=\{1\}$.
Then $a=2,c=1$. From Lemma \ref{main2}, on $O(6)$,
$$\sum_{\substack{i_1,i_2\in\cI,i_1,i_2\ne1\\i_1<i_2}}
\left[\begin{array}{cc}\ov1&\bi_2\\ i_1&\bi_1\\i_2\end{array}:
              \begin{array}{cc}1&2\\\ov2&3\\\ov3\end{array}\right]=
-\left[\begin{array}{cc}\ov1&\ov1\\ 1\end{array}\!:
              \begin{array}{cc}1&3\\\ov3\end{array}\right]
-\left[\begin{array}{cc}\ov1&\ov1\\ 1\end{array}:
              \begin{array}{cc}1&2\\\ov2\end{array}\right]+[\ov1:1].$$
The sum on the left is equal to
$$\left[\!\begin{array}{cc}\ov1&\!\ov2\\ \ov2&\!2\\2\end{array}\!\!:
              T\right]+
\left[\!\begin{array}{cc}\ov1&3\\ \ov2&2\\ \ov3\end{array}\!\!:T\right]+
\left[\!\begin{array}{cc}\ov1&\ov3\\ \ov2&2\\3\end{array}\!\!:T\right]+
\left[\!\begin{array}{cc}\ov1&3\\ 2&\ov2\\\ov3\end{array}\!\!:T\right]
+\left[\!\begin{array}{cc}\ov1&\ov3\\ 2&\ov2\\3\end{array}\!\!:T\right]
+\left[\!\begin{array}{cc}\ov1&\ov3\\ \ov3&3\\ 3\end{array}\!\!:T\right].$$
Then we get, on $O(6)$,
$$\left[\!\begin{array}{cc}\ov1&\ov2\\ \ov2&2\\2\end{array}\!\!:
              T\right]=
-\left[\!\begin{array}{cc}\ov1&2\\ \ov2&3\\ \ov3\end{array}\!\!:T\right]-
\left[\!\begin{array}{cc}\ov1&2\\ \ov2&\ov3\\3\end{array}\!\!:T\right]-
\left[\!\begin{array}{cc}\ov1&\ov2\\ 2&3\\\ov3\end{array}\!\!:T\right]
-\left[\!\begin{array}{cc}\ov1&\ov2\\ 2&\ov3\\3\end{array}\!\!:T\right]$$
$$\hfill
-\left[\!\begin{array}{cc}\ov1&\ov3\\ \ov3&3\\ \ov3\end{array}\!\!:T\right]
-\left[\begin{array}{cc}\ov1&\ov1\\ 1\end{array}\!:
              \begin{array}{cc}1&3\\\ov3\end{array}\right]
-\left[\begin{array}{cc}\ov1&\ov1\\ 1\end{array}:
              \begin{array}{cc}1&2\\\ov2\end{array}\right]
+[\ov1:1].$$

Next we deal with violations of (OS 2)
\begin{lemma}\label{3.7}
Suppose that $S$ and $T$ are column-increasing $\lambda$-tableaux
with two columns, and
that $S$ is not $O(n)$-standard, in that $\alpha_j+\beta_j=2j$ for some $j$
with $\alpha_j>\beta_j$, and $S(\alpha_j,1)=j$,  
$S(\beta_j,2)=\ov\jmath$, $S(\alpha_j-1,1)\ne\bj$,
Then $[S:T]$ may be expressed as a signed sum of bideterminants 
 as in equation (\ref{equn3.6}), Lemma \ref{3.6}.
\end{lemma}
\begin{proof}
Define $\cJ,\cK,A,B,C,a,b,c$ as in the proof of Lemma \ref{3.6}.
Since $j\in S_1$ and $\bj\in S_2$, then $j\in A$.

We claim that $\ov\jmath\in C$. If $\ov\jmath$ is in $S_1$, 
since $S_{\alpha_j,1}=j$,
 we must have $\ov\jmath$ in position 
$\alpha_j-1$ of column 1, since the entries of the columns are strictly
increasing. But we are assuming that  $S(\alpha_j-1,1)\ne\bj$.
So $\bj$ is not in column 1 of $S$. 
Since $S(\beta_j,1)=\bj$, and $\bj<j$, 
then $j$ does not occur in $S_2$.
So $\bj\in C$, as claimed.

Since $\alpha_j+\beta_j=2j$, then $2a+b=2j$, and since $a+b+c=2j$,
then $a=c$. Let $A'=A\cup\{\bj\}$, $C'=C-\{\bj\}$, and let $a'=\#A',
c'=\#C'$, so $c'<a$.  

Now proceed with the argument of the proof of Lemma 5.1, but use $C'$ instead of
$C$ in the definition of $\cS$. More precisely,
 let $\cS'$ be the set of all tableaux obtained from $S$ by
 replacing the pairs $i,\bi$ in $S$ by $i',\bi'$ where $i'\in \cI-C'$.

From Lemma \ref{main2}, the sum $\sum_{U\in\cS'}[U:T]$
is equal to a signed sum of bideterminants having shapes $\mu$ 
with $|\mu|<|\lambda|$. One of the terms in the sum $\sum_{U\in\cS'}[U:T]$
 is the original $[S:T]$. Of the others, all but one involve the 
replacement of at least one pair $i,\bi$ with $i',\bi'$ where $i'\in \cK$,
and so
the resulting $U\succ S$. The one exception arises from
the replacement
of $j,\bj$ with $\bj,j$. The resulting $U\succ S$ since
the second column of $U$ is identical to that of $S$ except that
$\bj$ has been replaced by $j$, and $j>\bj$. This completes the proof.
\end{proof}
\bigskip
\noindent Example. Consider the $O(7)$ bideterminant
$$[S:T]=\left[\begin{array}{cc}\ov1&\ov2\\1\\2\end{array}:
              \begin{array}{cc}1&2\\\ov2\\2\end{array}\right].$$
For $S$, with $j=2$, $\alpha_2+\beta_2=3+1=2j$ and $\alpha_2>\beta_2$;
$S(\alpha_2,1)=2, S(\beta_2,2)=\ov2$, but the 2 in the first column of $S$
is not protected by a $\ov2$ above it. 
(Note that for $T$, the 2 in the first column
is protected, and $T$ is $O(7)$-standard.)

Using the notation of the lemma, applied to $S$, with $j=2$,
we have $A=\{2\}$, $C=\{\ov2\}$, $C'$ is empty.
Let
$$s=\sum_{i\in\{\ov1,1,\ov2,2,\ov3,3,0\}}
\left[\begin{array}{cc}\ov1&i\\1\\\bi\end{array}:
              \begin{array}{cc}1&2\\\ov2\\2\end{array}\right].
$$
From Lemma \ref{main2}, as functions on $O(7)$,
$$s=\left[\begin{array}{c}\ov1\\1\end{array}:
              \begin{array}{c}1\\2\end{array}\right].$$
The sum $s$ is equal to
$$\left[\begin{array}{cc}\ov1&2\\1\\\ov2\end{array}:T\right]+
              \left[\begin{array}{cc}\ov1&\ov2\\1\\2\end{array}:T\right]+
\left[\begin{array}{cc}\ov1&3\\1\\\ov3\end{array}:T\right]+
\left[\begin{array}{cc}\ov1&\ov3\\1\\3\end{array}:T\right]
+\left[\begin{array}{cc}\ov1&0\\1\\0\end{array}:T\right].
$$
It follows that on $O(7)$,
$$\left[\begin{array}{cc}\ov1&\ov2\\1\\2\end{array}:T\right]=
-\left[\begin{array}{cc}\ov1&2\\1\\\ov2\end{array}:T\right]
-\left[\begin{array}{cc}\ov1&3\\1\\\ov3\end{array}:T\right]
-\left[\begin{array}{cc}\ov1&\ov3\\1\\3\end{array}:T\right]
-\left[\begin{array}{cc}\ov1&0\\1\\0\end{array}:T\right]
+\left[\begin{array}{c}\ov1\\1\end{array}:
              \begin{array}{c}1\\2\end{array}\right].$$

We now deal with violations of condition (OS 3). 

\begin{lemma}\label{3.8} 
Suppose that $S$ and $T$ are $\lambda$-tableaux with two columns. and
that $S$ is not $O(n)$-standard, in that $\alpha_j+\beta_j=2j$ for some $j$
with $\alpha_j=\beta_j$, and $S(\alpha_j,1)=\bj$, $S(\alpha_j,2)=j$ 
and $S(\alpha_j-1,2)\ne\bj$. Then $[S:T]$ can be expressed as
one-half a signed sum of bideterminants of three types:
(i) of shape $\mu$ where $|\mu|<|\lambda|$;
(ii) of shape $\mu$ where $\mu\lhd\lambda$;
(iii) of the form $[U:T]$ where $S\prec U$.
The tableaux $U$ in (iii) and their signs in the signed sum are independent
of $T$. The terms from (i) and (ii) are all 0 if $T=T^\lambda$.
\end{lemma}
\begin{proof}
Proceed as in the proof of the previous lemma.
Since $S(\alpha_j,1)=\bj$ and $j>\bj$,
then $j$ does not occur in the first column of $S$; 
since $S(\alpha_{j},2)=j$ and
$S(\alpha_j-1,2)\ne\bj$, 
then $\bj$ does not occur in the second column of $S$.
So $j\in C$; in this case let $C''=C-\{j\}$. Argue as in the proof of
the previous lemma, but now with $\cS''$ given by replacing pairs $i,\bi$ in
$S$ by $i',\bi'$ where $i'\in\cI-C''$. This time,
 $\sum_{U\in\cS''}[U:T]$ is equal to $[S:T]$ plus an exceptional term
$[S^*:T]$ where $S^*$ comes from replacing $\bj,j$ in $S$ with $j,\bj$,
plus a signed sum $s_3$ of bideterminants $[U:T]$ where $S\prec U$. From
Lemma \ref{main2}
$$[S:T]+[S^*:T]+s_3=s_1$$
where $s_1$ is a signed sum of bideterminants of shape $\mu$ 
where $|\mu|<|\lambda|$.
The tableau $S^*$ is not $GL(n)$-standard; straightening $[S^*:T]$, from
Lemma \ref{one-switch}, gives
$$[S^*:T]=[S:T]+s'_3+s_2$$
where $s'_3$ is a signed sum of bideterminants
$[V:T]$ of shape $\lambda$ where $V\succ S$, 
and $s_2$ is a signed sum of bideterminants
of shape $\mu$ where $\mu\lhd\lambda$.
Then
$$2[S:T]+s_3+s_3'+s_2=s_1,\quad 
2[S:T]=s_1-s_2-(s_3+s'_3),$$
where $s_1$ satisfies (i) in the conclusion of the lemma, $s_2$
satisfies (ii), and $s_3+s'_3$ satisfies (iii).  
\end{proof}
\noindent Example.
Consider the function on $O(6)$ given by
$$[S:T]=\left[\begin{array}{cc}1&1\\ \ov 2& 2\\ 3\end{array}:
              \begin{array}{cc}\ov1&   1\\ \ov2 &2 \\ 2\end{array}\right].$$
For the tableau $S$, $\alpha_2=\beta_2=2$. The $2$ in $S$ is
not protected, since $\ov2, 2$ occur in the second row, and there is no $\ov2$
above the 2. From Lemma \ref{main2}, on $O(6)$
$$\sum_{i\in\cI}
\left[\begin{array}{cc}1&1\\ i& \ov i\\ 3\end{array}:
              \begin{array}{cc}\ov1&   1\\ \ov2 &2 \\ 3\end{array}\right]
=-\left[\begin{array}{cc}1&1\\ 3\end{array}:
              \begin{array}{cc}\ov1&   1\\ 3\end{array}\right]
-\left[\begin{array}{cc}1&1\\ 3\end{array}:
              \begin{array}{cc}\ov2&   2\\ 3\end{array}\right].$$
Let the right side of this equation be $s_1$.
The left side of the equation is equal to
$$\left[\begin{array}{cc}1&1\\ \ov2& 2\\ 3\end{array}:
              \begin{array}{cc}\ov1&   1\\ \ov2 &2 \\ 3\end{array}\right]
+\left[\begin{array}{cc}1&1\\ 2& \ov 2\\ 3\end{array}:
              \begin{array}{cc}\ov1&   1\\ \ov2 &2 \\ 3\end{array}\right]
+\left[\begin{array}{cc}1&1\\ \ov3& 3\\ 3\end{array}:
              \begin{array}{cc}\ov1&   1\\ \ov2 &2 \\ 3\end{array}\right].$$
The first term is $[S:T]$. Call the third term $[V:T]$.  
The middle term needs to be $GL(n)$-straightened;
this was done in Example \ref{mead-eg}. From this example,
the middle term is equal to $[S:T]$, minus a bideterminant
$[U:T]$, where $S\prec U$, plus a signed sum $s_2$ of bideterminants
of shape $(2,1^3)$.
So we have, on $O(6)$,
$$2[S:T]=-[U:T]+[V:T]+s_1+s_2=s_1+s_2+s_3\mbox{ where }s_3=-[U:T]+[V:T].$$

We can now state our main straightening result.
\begin{theorem}\label{O(n)-straight}
Suppose that $S$ and $T$
are tableaux of the same shape $\lambda$, having at most $n$ parts. Then
as functions on $O(n)$,
$$[S:T]=\sum_Ua_U[U:T]+s$$
where the tableaux $U$ in the sum are $O(n)$-standard, each $a_U\in K$
and is independent of $T$, and $s$ is a linear combination of bideterminants
of shapes $\mu$ where $|\mu|<|\lambda|$ or $|\mu|=|\lambda|$ and 
$\mu\lhd\lambda$. If $T=T^\lambda$ then $s=0$.
\end{theorem}
\begin{proof} We use triple induction, first on the size $|\lambda|$
of a partition $\lambda$, next using the dominance order $\lhd$
on partitions of the same size, and lastly, downward induction on the
order $\prec$ on tableaux of the same shape. Given a bideterminant
$[S:T]$, use Mead's $GL(n)$-straightening from Theorem \ref{GL(n)-straight},
to reach the case that $[S:T]$ is $GL(n)$-standard. 

To achieve $O(n)$-standardness, we write the bideterminant $[S:T]$
as a product
of two bideterminants, the first $[S:T]_1$ coming from the first two columns
of each of $S$ and $T$, except for violations of (OS3) where we take
columns 1 and $b$, and the second $[S:T]_2$
coming from the remaining columns.
We need to straighten $[S:T]_1$.
To get the
condition that $\lambda'_1+\lambda'_2\le n$, use Lemma \ref{tilde-lambda}.
To achieve the three conditions in the definition of $O(n)$-standard,
use the three previous lemmas. Then  multiply by $[S:T]_2$.
In the course of using these results,
we sometimes replace a two-column subtableau with one where the columns
have different lengths. We might get an array which is not a tableau,
since the column lengths might not be decreasing. Then we rearrange
the columns to get a tableau.
\end{proof}

Using Remark \ref{transpose}, we have an analogous result where we replace
$[S:T]$ by $[T:S]$. Then we have
\begin{cor}The coordinate ring $K[O(n)]$ is spanned by $O(n)$-standard
bideterminants.
\end{cor}
\section{Filtrations and Linear Independence}

Let $K$ be a field of characteristic 0.

The coordinate ring $K[GL(n)]$ is a $GL(n)$-bimodule, 
$$g_1fg_2(x)=f(g_2xg_1),\quad f\in K[GL(n)],~ g_1,g_2,x\in GL(n).$$
The same equation defines an $O(n)$-bimodule structure on the 
coordinate ring $K[O(n)]$.

\begin{lemma}\label{binet-cauchy}
For tableaux $S$ and $T$ of the same shape $\lambda$, and $g\in GL(n)$,
$$[S:T]g=\sum_Ua_U[U:T]$$
where the sum is over a certain collection $\{U\}$ of $\lambda$-tableaux,
and each $a_U\in K$ does not depend on $T$.
\end{lemma}
\begin{proof}
For an $n\times n$
 matrix $A$, let $[S:T](A)$ denote the evaluation of the function $[S:T]$
at $A$; in this notation, $([S:T]g)(X)=[S:T](gX)$.

First supppose that $S$ and $T$ each have just one column. 
From the Binet-Cauchy formula, (see for example [Pra, 2.3, p.~10]])
$$[S:T](gX)=\sum_U[S:U](g)[U:T](X),\quad\mbox{so~}
[S:T]g=\sum_U\Bigl([S:U](g)\Bigr)[U:T]$$
where the $U$ in the sum vary over all column-increasing tableaux, 
of the same shape as $S$. 
In the general case, since a bideterminant is a product of determinants,
one for each column of the tableaux, $[S:T]=\prod_j[S_j:T_j]$.
Then
\begin{eqnarray*}[S:T]g&=&\prod_j[S_j:T_j]g=
\prod_j\left(\sum_{U_j}\Bigl([S_j:U_j](g)\Bigr)[U_j:T_j]\right)\\
&=&\sum_U\left(\prod_j \Bigl([S_j:U_j](g)\Bigr)[U_j:T_j]\right)\\
&=&\sum_U\left(\prod_j [S_j:U_j](g)\right)
\left(\prod_j [U_j:T_j]\right)\\
&=&\sum_Ua_U[U:T]
\end{eqnarray*}
where
$$a_U=\prod_j[S_j:U_j]g=[S:U]g.$$
\end{proof}

For a partition $\lambda$ having at most $n$ parts,
the  {\em Schur module} for $GL(n)$, which we denote by $L^\lambda$, is 
the $K$-span of all bideterminants $[T^\lambda:T]$
where $T$ varies over all tableaux of shape $\lambda$.
As is well known, this is a left $GL(n)$-submodule of $K[GL(n)]$;
this also follows from Lemma \ref{binet-cauchy}.
A good reference for Schur modules is Chapter 4 of [G], where
they are denoted by $D_{\lambda,K}$.  Since 
the characteristic of $K$ is 0, $L^\lambda$  is irreducible.
Define the right Schur module ${}^\lambda L$ to be the $K$-span
of all bideterminants $[T:T^\lambda]$.

Define the 
{\em orthogonal} left Schur module, denoted $L_O^\lambda$
to be the $K$-span of $[T^\lambda:T]_O$ 
as $T$ varies over all tableaux of shape 
$\lambda$. This is a left $O(n)$-submodule of $K[O(n)]$. Similarly we
have the right orthogonal
Schur module $[T:T^\lambda]_O$, spanned by  $[T^\lambda:T]_O$.

\begin{lemma}\label{schur}The right orthogonal
 Schur module ${}^\lambda L_O$ is spanned over $K$ by
all $O(n)$-standard bideterminants $[T:T^\lambda]$. The left
analogue $L^\lambda_O$ is spanned by all $O(n)$-standard 
bideterminants $[T^\lambda:T]$.
\end{lemma}
\begin{proof} The first statement follows from Theorem \ref{O(n)-straight}.
The second statement follows by taking transposes.
\end{proof}

We define the module $M^\lambda$, as in [KW, p.~254], as follows.
Fix a partition $\lambda$, and let $T$ be a $\lambda$-tableau.
Identify $T$ with an element of tensor space $V^{\otimes |\lambda|}$,
and let $\{T\}=Y^\lambda T$ where $Y^\lambda$ is the Young symmetrizer,
as in [KW, (2.2), p.~254]. Let $M^\lambda$ be the span of all $Y^\lambda T$,
as $T$ varies over all $\lambda$-tableaux.
This is a module for $GL(n)$, which as is well-known, is irreducible.

A {\em trace tensor} in $V^\ell$ is a linear combination of the form
$$\sum_{i\in\cI}x\otimes v_i\otimes y\otimes v_{\bi}\otimes z$$
where $x$, $y$ and $z$ are elements of some (possibly zero) tensor power
of $V$ and $\{v_i:i\in\cI\}$ is the standard basis of $V$.
 Let $U\subset V^{\otimes\ell}$ be the span of all such trace tensors.
Define
$$O^\lambda=M^\lambda/(M^\lambda\cap U).$$
As in [KW] and [W, Chap V.7], $O^\lambda$ is an irreducible $O(n)$-module, and
the irreducible $O(n)$-modules
are given by $O^\lambda$ where $\lambda$ is a partition with
at most $n$ parts, with $\lambda_1'+\lambda_2'\le n$.

In the $GL(n)$ case, it is well known that 
$M^\lambda$ and $L^\lambda$ are isomorphic. 
We will show that
$O^\lambda$ and $L_O^\lambda$ are isomorphic $O(n)$-modules. 
To do this we use the following formalism.
Let $M$ be a finite-dimensional polynomial $O(n)$-module, with basis
$m_1,\ldots,m_k$. For $m\in M$, $g\in O(n)$, write
$$gm=\sum_{i=1}^ka_i(g,m)m_i.$$
It is clear that for $g,g'\in O(n)$, we have $a_i(gg',m)=a_i(g,g'm)$.
Define $f_m\in K[O(n)]$ by $f_m(g)=a_1(g,m)$. (This, of course, depends on
the choice of ordered basis of $M$.) Define the map
$$F:M\to K[O(n)]\quad F(m)=f_m, \quad m\in M.$$
Then $F$ is a homomorphism of $O(n)$-modules.
We sometimes write $f_{M,m}$ instead of just $f_m$, and $F_M$ instead of
$F$.

If we have another $O(n)$-module $M'$ with basis 
$m_1'.\ldots,m_\ell'$, we use the
ordered basis $\{m_i\otimes m'_j\}$ whose first element is $m_1\otimes m'_1$.
For the element $m\otimes m'$, we have
$$g(m\otimes m')=gm\otimes gm'=\sum_{i,j}a_i(g,m)a'(g,m')m_i\otimes m'_j$$
from which it follows that
$$f_{M\otimes M',m\otimes m'}=f_{M,m}f_{M',m'}.$$

If $N$ is an $O(n)$-sumbodule of $M$ which does not contain a scalar multiple
of $m_1$, 
 we pick the basis $m_1,\ldots,m_k$
of $M$ such that for some $j>1$, $m_{j+1},\ldots,m_k$ is a basis of $N$.
Let $\ov M$ denote $M/N$, which has basis $\ov m_1,\ldots, \ov m_j$.
If $m\in M$ has the form $m=\sum_{i=1}^ka_im_i$ then $\ov m=
\sum_{i=1}^ja_i\ov m_i$ so $f_{\ov M,\ov m}=f_{M,m}$.

Suppose that $\lambda$ is the partition of $k$ whose shape consists of one
column.
It is well-known and easy to see that the module $GL(n)$-module 
$M^\lambda$ is the $k$-th
exterior power $\bigwedge^k V$ of $V$; $M^\lambda$ is an $O(n)$-module
by restriction.
Pick a basis of $M^\lambda$ whose
first element is $m_1=v_{j_1}\wedge\cdots\wedge v_{j_k}$ where $j_1,j_2,\ldots
j_k$ are the first $k$ elements of $\cI$, as in (\ref{order}).  Let
$T(k)$ denote the one-column tableau with entries $j_1,j_2,\ldots j_k$. Then
$f_{M^\lambda,m_1}=[T(k):T(k)]_O$.  Let $\wedge v(k)=v_{j_1}\wedge\cdots\wedge
v_{j_k}$.

Now let $\lambda$ be any partition. Then $M^\lambda$ 
has a basis whose first element
$m_1$ is $Y^\lambda T^\lambda$ which is equal to
$$\wedge v(\lambda'_1)\otimes\cdots\otimes\wedge v(\lambda'_\ell)
\in 
\bigwedge\nolimits^{|\lambda'_1|}V\otimes\cdots\otimes 
\bigwedge\nolimits^{|\lambda'_\ell|}V$$
where $\lambda'$ has $\ell$ parts.
The module $M^\lambda$ is an $O(n)$-module by restriction. 
It follows that $f_{M^\lambda,m_1}$ is the bideterminant 
$[T^\lambda:T^\lambda]_O$.

Since $O^\lambda$ is a factor module of $M^\lambda$, it follows
that $f_{O^\lambda,\ov m_1}=f_{M^\lambda,m_1}=[T^\lambda:T^\lambda]_O$.
Since $O^\lambda$ is an irreducible $O(n)$-module, then any
element $m\in O^\lambda$ is a linear combination of elements of the
form $\{g_i\ov m_1\}$ for suitable $g_i\in O(n)$. Hence
$F_{O^\lambda}(m)\in L_O^\lambda$.
So $F_{O^\lambda}$ is an $O(n)$-module homomorphism 
from $O^\lambda$ to $L_O^\lambda$. From [KW], $O^\lambda$ has basis given by
standard $O(n)$-tableaux. Note that the main idea of [KW] is to prove
straightening; linear independence follows from Proctor [Pro].

\begin{prop} \label{basis1}For a partition $\lambda$ with at most $n$ parts,
the left and right Schur modules 
$L_O^\lambda$  and ${}^\lambda L_O$
are irreducible, and have $K$-bases 
given by $O(n)$-standard bideterminants $[T^\lambda:T]$, and $[T:T^\lambda]$,
respectively. 
\end{prop}
\begin{proof} From the preceding paragraph there is a non-zero homomorphism
$F$ from $O^\lambda$ to $L_O^\lambda$. 
Since $O^\lambda$ has $K$-basis indexed by the 
set of $O(n)$-standard tableaux of shape 
$\lambda$,
and since $L_O^\lambda$ is spanned by $O(n)$-standard 
$[T^\lambda:T]_O$  
of shape $\lambda$, by Lemma \ref{schur},
it follows that $F$ is surjective. It is injective, since $O^\lambda$
is irreducible. The analogue for right Schur modules is similar.
\end{proof}

Let $\le$ be any partial order on the set of partitions of at most
$n$ parts, such that 
\begin{equation}\label{<}
\mu<\lambda \mbox{ if }|\mu|<|\lambda|;\mbox{ and }
\mu<\lambda\mbox{ if }|\mu|=|\lambda|\mbox{ and }\mu\lhd\lambda.
\end{equation}

Let $A(\le \lambda)$ denote the $K$-subspace of $K[O(n)]$ spanned by
standard $O(n)$-bideterminants of shape $\mu$ where $\mu\le\lambda$.
Define $A(<\lambda)$ similarly.
Let $\ov{A(\le\lambda)}$  denote $A(\le \lambda)/A(< \lambda)$.
Let $\ov{[S:T]_O}$ denote the image of the bideterminant $[S:T]_O$ in 
$A(\le \lambda)/A(< \lambda)$.
\begin{theorem}\label{hom} There is an isomorphism 
$$\Phi:L^\lambda_O\otimes {}^\lambda L_O\to \ov{A(\le\lambda))}$$
of $O(n)$-bimodules given by
$$\Phi\Bigl([T^\lambda:S]_O\otimes [T:T^\lambda]_O\Bigr)=\ov{[T:S]_O}$$
where $S$ and $T$ are $O(n)$-standard tableaux of shape $\lambda$.
\end{theorem}
\begin{proof}
Note that $\Phi$ is well-defined, since $L^\lambda_O$ and $ {}^\lambda L_O$
have bases, from Proposition \ref{basis1}.
We first show that $\Phi$ is a right $O(n)$-module homomorphism. 

Suppose that $S$ and $T$ are $O(n)$-standard.
For $g\in O(n)$, we have from Lemma \ref{binet-cauchy}
\begin{equation}\label{a_U} [T:T^\lambda]g=\sum_Ua_U[U:T^\lambda]
\end{equation}
where the $U$ in the sum are column-increasing tableaux, and $a_U\in K$.
From Theorem \ref{O(n)-straight}, each $[U:T^\lambda]_O$ can be
written as
\begin{equation}\label{b} [U:T^\lambda]_O=\sum_{V}b_{U,V}[V:T^\lambda]_O
\end{equation}
where the tableaux $V$ are $O(n)$-standard.  Then
\begin{eqnarray}\label{c}
\Phi\Bigl([T^\lambda:S]_O\otimes [T:T^\lambda]_Og\Bigr)&=&
\Phi\Bigl([T^\lambda:S]_O\otimes \sum_{U,V}a_Ub_{U,V}[V:T^\lambda]_O\Bigr)
\nonumber\\
&=&
\sum_{U,V}a_Ub_{U,V}\ov{[V:S]_O}.
\end{eqnarray}
On the other hand,
$$\Bigl(\Phi\left([T^\lambda:S]_O\otimes [T:T^\lambda]_O\right)\Bigr)g=
\ov{[T:S]_Og}.$$
From Lemma \ref{binet-cauchy}
the coefficients $a_U$ in (\ref{a_U}) are the same for
$[T:S]g$ as they are for $[T:T^\lambda]g$, so
$$[T:S]g=\sum_Ua_U[U:S].$$
Moreover when $O(n)$-straightening $[U:S]_O$, we get the same answer, 
mod $A(<\lambda)$,
as we do when straightening $[U:T^\lambda]_O$, from Theorem 
\ref{O(n)-straight}. So from (\ref b) we get
$$[U:S]_O\equiv \sum_Vb_{U,V}[V:S]_O \mbox{ mod }A(<\lambda).$$
Then
$$\Bigl(\Phi\left([T^\lambda:S]_O\otimes [T:T^\lambda]_O\right)\Bigr)g=
\ov{[T:S]_Og}=\sum_{U,V}a_Ub_{U,V}\ov{[V:S]_O}.$$
It follows from this and (\ref c) that $\Phi$ is a homomorphism
of right $O(n)$-modules. That it is a homomorphism of
left $O(n)$-modules has a similar proof.

Since $L^\lambda_O$ is an irreducible
left $O(n)$-module and  ${}^\lambda L_O$ is an irreducible
right $O(n)$-module, then $L^\lambda_O\otimes {}^\lambda L_O$
is an irreducible $O(n)$-bimodule. Therefore $\Phi$ has trivial kernel.
Now $L^\lambda_O\otimes {}^\lambda L_O$ has $K$-basis given by the set of
all $[T^\lambda:S]\otimes [T:T^\lambda]$
where $S$ and $T$ are $O(n)$-standard of shape $\lambda$,
and this basis is mapped
by $\Phi$ to the set $\{[T:S]_O\}$ of generators of $\ov{A{(\le \lambda)}}$. 
It follows that $\Phi$ is an isomorphism. 
\end{proof}

\begin{cor}\label{hom-cor}The set of $O(n)$-standard
bideterminants $[S:T]_O$ are linearly independent.
\end{cor}
\begin{proof}
This is because the set of all $[T^\lambda:S]_O\otimes[T:T^\lambda]$
is a basis of $L^\lambda_O\otimes {}^\lambda L_O$.
\end{proof}

\begin{theorem}\label{lin-ind}
Let $<$ be a partial order satisfying (\ref<), and let $\tilde <$ be a total
order which refines $<$. Then $A(\tilde\le\lambda)$ has $K$-basis given
by the set of $O(n)$-standard bideterminants $[S:T]_O$ of
shape $\mu\tilde\le\lambda$. 
\end{theorem}
\begin{proof}
By induction, $A(\tilde<\lambda)$
has a basis given by $O(n)$-standard bideterminants of shape $\tilde<\lambda$.
It follows from Theorem \ref{hom} that
the images of $[S:T]$ in $A(\le\lambda)/A(<\lambda)$, 
where $[S:T]$
are $O(n)$-standard of shape $\lambda$, are a basis of    
$A(\le\lambda)/A(<\lambda)$. 
We conclude that the $O(n)$-standard bideterminants
of shape $\tilde\le\lambda$ form a basis of $A(\tilde\le\lambda)$. 
\end{proof}
\begin{cor}
The set of $O(n)$-standard bideterminants
is a $K$-basis of $K[O(n)]$.
\end{cor}
\begin{proof}
From Theorem \ref{O(n)-straight} the  $O(n)$-standard bideterminants span
$K[O(n)]$. Finitely many of them lie in $A(\tilde\le\lambda)$ 
for some $\lambda$,
hence are linearly independent.
\end{proof}

\section{Base Change and Non-zero Characteristic}

We first discuss base change. 
Suppose that $R$ is an integral domain of characteristic not 2,
having infinite field of fractions $K$.
Define $O(n,R)$ to be the subgroup of $GL(n,R)$ given by elements
$g$ such that $g^tJg=J$, where $J$ is the matrix of the bilinear form
(\ref{bilform}). 
Let $SO(n,R)$ be the elements in $O(n,R)$ of determinant 1. Essentially
following [Bo1, 2.4] define $R[O(n)]$ to be the polynomial ring
$R[X(i,j)]$ modulo the ideal $\mathcal A$ of polynomials which vanish on
$O(n,\ov K)$ where $\ov K$ is the algebraic closure of $K$. Define
$R[SO(n)]$ similarly. 
From [Bo2, 18.3], $SO(n,K)$ is dense in $SO(n,\ov K)$, in the Zariski
topology, and then $O(n,K)$ is dense in $O(n,\ov K)$. So $\mathcal A$
is the ideal of polynomials in  $R[X(i,j)]$ which vanish on $O(n,K)$. 
In the case that $R$ is an infinite field $K$, the definition of $K[O(n)]$
just given agrees with the definition given  
in section 2, namely the restriction of polynomials in $K[X(i,j)]$ to
$O(n,K)$.

The group $SO(n,K)$ is a Chevalley group. (Strictly speaking, for this
to hold we should assume that $n\ge 4$.)  We use a result of
Chevalley, as discussed in [Bo1]. 
From [Bo1, Lemma 4.6] the ring 
$K\otimes_{\Z}\Z[SO(n)]$ is reduced (that is, it has no nilpotent
elements). It follows, as in [Bo1, 3.4], that
\begin{equation}\label{SO-base-change}
K[SO(n)]\cong K\otimes_{\Z}\Z[SO(n)].
\end{equation}
The group $O(n,K)$ has two irreducible components, namely the
subsets on which the determinant is $\pm 1$. 
We have the idempotents $e_+=(1/2)(1-\det)$ and 
 $e_-=(1/2)(1+\det)$ in $K[O(n)]$, and
$$K[O(n)]=e_+K[O(n)]\oplus  e_-K[O(n)],\qquad e_+K[O(n)]=K[SO(n)].$$
Let $\Z[1/2]$ be the ring of rational numbers whose denominators are powers
of 2.
It follows that
\begin{equation}\label{base-change}
K[O(n)]\cong K\otimes_{\Z[1/2]}\Z[1/2][O(n)].
\end{equation}

%

\begin{theorem}\label{any-R}If $K$ is an infinite field of odd characteristic,
then the set of $O(n)$-standard bideterminants $[S:T]$
 is a $K$-basis of $K[O(n)]$.
\end{theorem}
\begin{proof}
We first prove an analogous result for $\Z[1/2][O(n)]$. 
Let $R=\Z[1/2]$.

Start with a polynomial $P$ in $R[X(i,j])$.
Mead straightening can be done over any field, indeed over $\Z$.
So $P$ is an $R$-linear combination
of bideterminants $[S:T]$ where $S$ and $T$ are $GL(n)$-standard.
A bideterminant $[S:T]$
is a polynomial in the $n^2$ variables $X(i,j)$ with integer coefficients,
 and so gives rise to an
an element, denoted, $[S:T]_{O,R}$,
 of $R[O(n)]$.
The $O(n)$-straightening results we have used, namely Lemmas 
\ref{tilde-lambda}, \ref{3.6}, \ref{3.7}, \ref{3.8} hold in $R[O(n)]$.
So we can do $O(n)$-straightening in $R[O(n)]$ and 
we see that the $O(n)$-standard bideterminants $[S:T]_{O,R}$
generate $R[O(n)]$ as an $R$-module.
In $\Q[O(n)]$,
the $O(n)$-standard bideterminants $[S:T]$ 
are linearly independent over $\Q$, from Theorem \ref{lin-ind},
 so the $O(n)$-standard bideterminants $[S:T]_{O,R}$
are linearly independent over $R$, from (\ref{base-change}). 
Thus $R[O(n)]$ is a free $R$-module,
with basis consisting of the $O(n)$-standard bideterminants $[S:T]_{O,R}$.
Using (\ref{base-change}), 
it follows that  the $O(n)$-standard bideterminants 
 $[S:T]$
form a $K$-basis  of $K[O(n)]$.
\end{proof}

We now consider Theorem \ref{hom} at odd positive characteristic.
The definitions of $L_O^\lambda$, ${}^\lambda L_O$, and 
$\ov{A\le\lambda}$ make sense for any field $K$. At 
positive characteristic, $L^\lambda_{O}$ is not in
general an irreducible left $O(n,K)$-module, and we must find
a different argument for the proof.

Let $R=\Z[1/2]$, and let
$L^\lambda_{O,R}$ be the $R$-span of $[T^\lambda:T]_{O,R}\in R[O(n)]$
for all $\lambda$-tableaux $T$. This is a 
left $O(n,R)$ module, and has a basis of bideterminants $[T^\lambda:T]_{O,R}$
where $T$ is $O(n)$-standard of shape $\lambda$.
There are analogous definitions and results
for ${}^\lambda L_{O,R}$.
  
As in the definition of $A(\lambda)$ following equation (\ref<), 
let $A(\le\lambda)_R$ be the $R$-span of all $[S:T]_{O,R}$ where
$S$ and $T$ have shape $\mu\le\lambda$. Similarly we define
$A(<\lambda)_R$, and let $\ov{A(\le\lambda)}_R=A(\le\lambda)_R/
A(<\lambda)_R$. Then we have the following analogue of Theorem \ref{hom}.
\begin{theorem} If $K$ is an infinite field of characteristic
not 2, then
there is an isomorphism 
$$\Phi:L^\lambda_{O}\otimes {}^\lambda L_{O}\to \ov{A(\le\lambda))}_{O}$$
of $O(n,K)$-bimodules given by
$$\Phi\Bigl([T^\lambda:S]_{O}\otimes [T:T^\lambda]_{O}\Bigr)=
\ov{[T:S]_{O}}$$
where $S$ and $T$ are $O(n)$-standard tableaux of shape $\lambda$.
\end{theorem}
\begin{proof}
That $\Phi$ is a bimodule homomorphism is similar to the proof 
of Theorem \ref{hom}.
The isomorphism $\Phi$ above is defined over $R=\Z[1/2]$:
$$\Phi_R:L^\lambda_{O,R}\otimes {}^\lambda L_{O,R}\to\ov{A(\le\lambda))}_{O,R},
\quad \Phi_R\Bigl([T^\lambda:S]_{O,R}\otimes [T:T^\lambda]_{O,R}\Bigr)=
\ov{[T:S]_{O,R}}.$$
Since $\Phi$ is injective if $K=\Q$, then
$\Phi_R$ is injective. It then follows that
$\Phi_R$ is an isomorphism, and by base change we see
that $\Phi$ is an isomorphism.
\end{proof}

\section{The group of orthogonal similitudes}

Let $R$  be as in the previous section.
With a view to future applications to the orthogonal Schur algebra,
we now formulate our results with respect to the group $GO(n,R)$
of {\em orthogonal similitudes}.
This is defined as
all $g\in GL(n,R)$ such that for some unit $\gamma=\gamma(g)$
of $R$,
\begin{equation}\label{bil-form-gamma}
\left<gv,gw\right>=\gamma\left<v,w\right>\mbox{ for all }v,w\in R^n.
\end{equation}
Here $\left<~,~\right>$ is the form defined in equation (\ref{bilform}).

The element $\gamma$ can be defined, as a function on
 $GO(n,R)$, by
\begin{equation}\label{gamma}
\gamma=\sum_{i\in\mathcal I}X(i,1)X(\bi,\ov1).\end{equation}
The map $\gamma:GO(n,R)\to R^\times$
is a homomorphism.

Define $R[GO(n)]$ to be the polynomials in the $n^2+1$ variables
$X(i,j)$ and $1/\det$, with coefficients in $R$, 
modulo the ideal of those polynomials vanishing on $GO(n,K)$.

\begin{lemma}\label{Gbase-change}
$$K[GO(n)]\cong K\otimes_{\Z[1/2]}\Z[1/2][GO(n)].$$
\end{lemma}
\begin{proof}
There are two cases, depending on whether $n$ is even or odd,
since $GO(n)$ is connected if $n$ is odd but not if $n$ is even.
For $g\in GO(n,K)$ the condition (\ref{bil-form-gamma}) is equivalent
to $g^tJg=\gamma J$,
where $J$ is the matrix of the bilinear form (\ref{bilform});
then on $GO(n,K)$,
\begin{equation}\label{det-gamma}
\det g^2=\gamma(g)^n.\end{equation}
If $n$ is odd, then from [KMRT, 12.4], $GO(n,K)=SO(n,K)\cdot K^{\times}\cong
SO(n,K)\times K^{\times}$.
It follows from the definitions that
$$\Z[GO(n)]\cong \Z[SO(n)]\otimes \Z[t,t^{-1}]$$
where $t$ is an indeterminate.
Then the lemma follows from (\ref{SO-base-change}).

Now suppose that $n$ is even. 
Given $c\in K^\times$ let $\xi(c)$ be the diagonal matrix
$$\mbox{diag}(c,c,\ldots,c,1,1,\ldots 1)$$
 where the first $n/2$ diagonal entries are equal to $c$.
We have a map of algebraic varieties 
$$\Phi:O(n)\times K^\times\to GO(n), \quad \Phi(g,c)=g\xi(c).$$
For $g\in GO(n,K)$, 
$$\gamma\bigl(g\xi(\gamma(g)^{-1}\bigr)
=\gamma(g)\gamma(g)^{-1}
=1,\quad\mbox{so }\quad
g\xi(\gamma(g)^{-1})\in O(n).$$
Then $\Phi$ is invertible, with inverse given by
$$\Psi(g)=(g\xi(\gamma(g))^{-1} ,\gamma(g))$$
and $\Psi$ is a regular map because of (\ref{gamma}).
It follows that
$$K[GO(n)]\cong K[O(n)]\otimes K[t,t^{-1}].$$
and the lemma follows from (\ref{base-change}).
\end{proof}

Define
$A(n,r,R)$ to be the $R$-module of polynomials of degree $r$ in the $n^2$,
variables $X(i,j)$, with coefficients in $R$.
Define
$A_{GO}(n,r,R)$ to be $A(n,r,R)$ modulo the submodule of all polynomials
which vanish on $GO(n,K)$, so
$A_{GO}(n,r,K)$ is the restriction of the polynomials in $A(n,r,K)$ to
$GO(n,K)$.
We will find a basis for $A_{GO}(n,r,K)$.

The main technical Lemma \ref{main}, becomes, in this context:
\begin{lemma}\label{main-g}As functions on $GO(n,K)$, we have
$$\sum_{\substack{i_1,\ldots,i_a\\ \in\cI-C}}
\left[\begin{array}{cc}i_1&\bi_1\\
              i_2&\bi_2\\
               \vdots&\vdots\\
              i_a&\bi_a\\
\multicolumn{2}{c}{S_0}
               \end{array}:T\right]
=a!\sum_{d=1}^a\gamma^{d}\cS_d.$$
\end{lemma}
\begin{proof} 
Proceed as with the proof of Lemma \ref{main}, until equation (\ref{main1}).
As functions on $O(n,K)$, $\sum_{i\in\cI}X(i,j)X(\bi,k)=\delta_{j,\ov k}$,
whereas as functions on $GO(n,K)$, we have
$$\sum_{i\in\cI}X(i,j)X(\bi,k)=\delta_{j,\ov k}\gamma.$$
Then equation (\ref{main1}) gets replaced, as functions on $GO(n,R)$, by
\begin{equation*}
\sum_{i_s\in\cI-C}
X\bigl(i_s,\sigma(g_s)\bigr)X\bigl(\bi_s,\tau(h_s)\bigr)
=\delta_{\sigma(g_s),\ov{\tau(h_s)}}\gamma
-\sum_{i_s\in C}X\bigl(i_s,\sigma(g_s)\bigr)
X\bigl(\bi_s,\tau(h_s)\bigr).
\end{equation*}
Similarly, all occurences of $\delta_{\sigma(g_s),\ov{\tau(h_s)}}$
get replaced by $\delta_{\sigma(g_s),\ov{\tau(h_s)}}\gamma$,
and $\delta(s)$ gets replaced by  $\delta_{\sigma(g_s),\ov{\tau(h_s)}}\gamma$.
Equation (\ref{e-main2}) becomes, on $GO(n,K)$,
\begin{equation*}
\prod_{s=1}^a\left(\delta(s)-\sum_{i_s\in C}Z(s,i_s)\right)
=\sum_{D\in\cD}(-1)^{\ell(D)}\gamma^d\sum_{(i_1,\ldots,i_\ell)\in C^\ell}
\prod_{s\in D'}Z(s,i_s).
\end{equation*}
With these replacements, the result follows as in the proof of Lemma 
\ref{main}.
\end{proof}

Our main result now, for any infinite field $K$ of characteristic not 2, is
\begin{theorem}A basis of $A_{GO}(n,r,K)$ is given by the set of functions on 
$GO(n,K)$ 
$$\cB_r=\{\gamma^k[S:T]:k\in\Z,
\quad 0\le k\le r/2\}$$
where $[S:T]$ is $O(n)$-standard of shape $\lambda$ with 
$|\lambda|=r-2k$.
\end{theorem}
\begin{proof}
As in the previous section, we first find a basis for 
$A_{GO}(n,r,\Z[1/2])$. 
We start with $GL(n)$-straightening, which can be done over any ring $R$.
We need to show that as functions on $GO(n,Q)$, a bideterminant $[S:T]$
where $S$ and $T$ are tableaux of shape $\lambda$ with $|\lambda|=r$,
is a $\Z[1/2]$-linear combination of the elements in $\cB_r$. 
We call this $GO(n)$-straightening.
We need $GO(n)$-analogues of our $O(n)$-straightening lemmas.
From Lemma \ref{main-g} we get a $GO$ analogue of
Lemma \ref{main2}, which is that on $GO(n,\Q)$
\begin{equation}\label{g-main2}
\sum_{U\in\cS^*}[U:T]=\sum_{d=1}^a\gamma^d\cS_d.
\end{equation}
As in Lemma \ref{3.6},
$[S:T]=\sum_{U\in\cS} -[U:T] + s$; using (\ref{g-main2}),
we see that 
on $GO(n,Q)$
$s$ is a signed sum of terms of the form 
$\gamma^b[S':T']$ where $S'$ is a tableaux of shape $\sigma$ 
where $|\sigma|<|\lambda|$ and $b+|\sigma|=r$. 
Since $[S':T']$ can be $GO(n)$-straightened
by induction, then so can $[S:T]$.
Analogues of Lemmas \ref{3.7} and \ref{3.8}
are similar. 

We need a $GO$-analogue of Lemma \ref{tilde-lambda}.
Equation (\ref{eqn-tilde-lambda}) becomes, on $GO(n,\Q)$,
$$[S:T]=\pm\det{}^2\gamma^\ell\cdot[\bar S_1':\bar T_1'][\bar S_2':\bar T_2']$$
for a suitable integer $\ell$. 
If $[S:T]$ has shape $\lambda$ where $\lambda_1'+\lambda_2'>n$,
the analogue of Lemma \ref{tilde-lambda}, using
(\ref{det-gamma}), is that on $GO(n,\Q)$, 
$[S:T]=\pm\gamma^b[\tilde S:\tilde T]$,
for a suitable integer $b$, 
where $\tilde S$ and $\tilde T$ are as in Lemma \ref{tilde-lambda}.
This  proves $GO(n)$-straightening. 

Since a polynomial in $\Z[1/2][X(i,j)]$ which vanishes on $GO(n,\Q)$
also vanishes on $O(n,\Q)$, we 
have a ring homomorphism
$$F:A_{GO}(n,r,Z[1/2])\to \Z[1/2][O(n)].$$  
The set $F(\cB_r)$ 
is linearly independent, 
since $\gamma=1$ on $O(n,\Q)$ and the $O(n)$-standard bideterminants
are linearly independent in $\Q[O(n)]$. 
So $\cB_r$ is a linearly independent set in $A_{GO}(n,r,Z[1/2])$.
The theorem now follows from Lemma \ref{Gbase-change}.
\end{proof}

\end{document}